\documentclass{article}
\usepackage{hyperref}
\usepackage{amssymb}
\usepackage{amsthm}
\usepackage{eucal}
\usepackage[dvips]{graphicx}
\usepackage{multirow}
\usepackage{fancyhdr}

\newtheorem{proposition}{Proposition}[section]

\newtheorem{thm}[proposition]{Theorem}

\newtheorem{lem}[proposition]{Lemma}

\newtheorem{rem}[proposition]{Remark}

\begin{document}

\title{Wiener Process with Reflection in Non-Smooth Narrow Tubes}
\author{Konstantinos Spiliopoulos\\
Lefschetz Center for Dynamical Systems\\
Division of Applied Mathematics\\
Brown University, Providence, RI, 02912\\
kspiliop@dam.brown.edu}

\date{}

\maketitle

\begin{abstract}
Wiener process with instantaneous reflection in narrow tubes of
width $\epsilon\ll 1$ around axis $x$  is considered in this paper.
The tube is assumed to be (asymptotically) non-smooth in the
following sense. Let $V^{\epsilon}(x)$ be the volume of the
cross-section of the tube. We assume that
$\frac{1}{\epsilon}V^{\epsilon}(x)$ converges in an appropriate
sense to a non-smooth function as
$\epsilon\downarrow 0$. This limiting function can be composed by smooth
functions, step functions and also the Dirac delta distribution. Under this assumption we prove that the
$x$-component of the Wiener process converges weakly to a Markov
process that behaves like a standard diffusion process away from the
points of discontinuity and has to satisfy certain gluing conditions at the
points of discontinuity.

\end{abstract}


\vspace{1cm}

\framebox{\parbox[c]{11cm}{This is an electronic reprint of the
original article published by the
\href{http://www.emis.de/journals/EJP-ECP/_ejpecp/search9ead.html}{Electronic
Journal of Probability, Vol. 14, Paper no. 69, pp. 2011-2037}. This
reprint differs from the original in pagination and typographic
detail.}}

\newpage

\section{Introduction}

For each $x\in \mathbb{R}$ and $0<\epsilon<<1$, let
$D^{\epsilon}_{x}$ be a bounded interval in $\mathbb{R}$  that
contains $0$. To be more specific, let $D^{\epsilon}_{x}=[-V^{l,\epsilon}(x), V^{u,\epsilon}(x)]$, where $V^{l,\epsilon}(x), V^{u,\epsilon}(x)$ are sufficiently smooth, nonnegative functions, where at least one of the two is a strictly positive function. Consider the state space
$D^{\epsilon}=\{(x,y):x\in \mathbb{R},y\in D^{\epsilon}_{x}\}\subset
\mathbb{R}^{2}$. Assume that the boundary $\partial D^{\epsilon}$ of
$D^{\epsilon}$ is smooth enough and denote by
$\gamma^{\epsilon}(x,y)$ the inward unit normal to $\partial
D^{\epsilon}$. Assume that $\gamma^{\epsilon}(x,y)$ is not parallel
to the $x$-axis.

Denote by $V^{\epsilon}(x)=V^{l,\epsilon}(x)+V^{u,\epsilon}(x)$ the length of the cross-section
$D_{x}^{\epsilon}$ of the stripe. We assume that $D^{\epsilon}$ is a
narrow stripe for $0<\epsilon<<1$, i.e. $V^{\epsilon}(x)\downarrow
0$ as $\epsilon \downarrow 0$. In addition, we assume that
$\frac{1}{\epsilon}V^{\epsilon}(x)$ converges in an appropriate
sense to a non-smooth function, $V(x)$, as $\epsilon\downarrow 0$.
The limiting function can be composed for example by smooth
functions, step functions and also the Dirac delta distribution. Next, we state the
problem and we rigorously introduce the assumptions on
$V^{\epsilon}(x)$ and $V(x)$. At the end of this introduction we
formulate the main result.

Consider the Wiener process $(X^{\epsilon}_{t},Y^{\epsilon}_{t})$ in $D^{\epsilon}$ with instantaneous normal reflection on the boundary of $D^{\epsilon}$. Its trajectories can be described by the stochastic differential equations:
\begin{eqnarray}
X_{t}^{\epsilon}&=& x+ W_{t}^{1}+\int_{0}^{t}\gamma_{1}^{\epsilon}(X_{s}^{\epsilon},Y_{s}^{\epsilon})dL_{s}^{\epsilon}\nonumber\\
Y_{t}^{\epsilon}&=& y + W_{t}^{2}+\int_{0}^{t}\gamma_{2}^{\epsilon}(X_{s}^{\epsilon},Y_{s}^{\epsilon})dL_{s}^{\epsilon}.\label{StochasticProcessWithReflection1}
\end{eqnarray}
Here  $W_{t}^{1}$ and $W_{t}^{2}$ are independent Wiener processes in $\mathbb{R}$  and $(x,y)$ is a point inside $D^{\epsilon}$; $\gamma_{1}^{\epsilon}$ and $\gamma_{2}^{\epsilon}$ are both projections of  the unit inward normal vector to $\partial D^{\epsilon}$ on the axis $x$ and $y$ respectively. Furthermore, $L^{\epsilon}_{t}$ is the local time for the process $(X^{\epsilon}_{t},Y^{\epsilon}_{t})$ on $\partial D^{\epsilon}$, i.e. it is a continuous, non-decreasing process that increases only when $(X^{\epsilon}_{t},Y^{\epsilon}_{t}) \in \partial D^{\epsilon}$ such that the Lebesgue measure $\Lambda\{t>0:(X^{\epsilon}_{t},Y^{\epsilon}_{t}) \in \partial D^{\epsilon}\}=0$ (eg. see \cite{KS}).

Our goal is to study the weak convergence of the $x-$component of
the solution to  (\ref{StochasticProcessWithReflection1}) as
$\epsilon\downarrow 0$. The $y-$component clearly converges to $0$
as $\epsilon\downarrow 0$. The problem for narrow stripes with a
smooth boundary was considered in \cite{F2} and in \cite{FS1}.
There, the authors consider the case
$\frac{1}{\epsilon}V^{\epsilon}(x)=V(x)$, where $V(x)$ is a smooth
function. It is proven that $X^{\epsilon}_{t}$ converges to a standard diffusion process $X_{t}$, as $\epsilon\downarrow 0$. More precisely, it is shown that for any $T>0$
\begin{equation}
\sup_{0\leq t \leq T}E_{x}|X_{t}^{\epsilon}-X_{t}|^{2}\rightarrow 0 \hspace{0.2cm} \textrm{as} \hspace{0.2cm} \epsilon \rightarrow 0,
\end{equation}
where  $X_{t}$ is the solution of the stochastic differential equation
\begin{equation}
X_{t}= x + W_{t}^{1} + \int_{0}^{t}\frac{1}{2}\frac{V_{x}(X_{s})}{V(X_{s})}ds \label{LimitingStochasticProcess1}
\end{equation}
and $V_{x}(x)=\frac{d V(x)}{dx}$.

In this paper we assume that
$\frac{1}{\epsilon}V^{\epsilon}(x)$ converges to a non-smooth
function as described by
(\ref{DefinitionOfCrossSections1})-(\ref{DefinitionOfCrossSections5})
below, as $\epsilon \downarrow 0$. Owing to the non smoothness of the limiting function, one cannot hope to obtain a limit in mean square sense to a standard diffusion process as before. In particular, as we will see, the non smoothness of the limiting function leads to the effect that the limiting diffusion may have points where the scale function is not differentiable (skew diffusion) and also points with positive speed measure (points with delay).

For any $\epsilon>0$, we introduce the functions
\begin{equation}
u^{\epsilon}(x):=\int_{0}^{x}2\frac{\epsilon}{V^{\epsilon}(y)}dy \hspace{0.3cm}\textrm{ and }\hspace{0.3cm}
v^{\epsilon}(x):=\int_{0}^{x}\frac{V^{\epsilon}(y)}{\epsilon}dy.\label{uANDvForEscortProcessIntro}
\end{equation}

Now, we are in position to describe the limiting behavior of $\frac{1}{\epsilon}V^{\epsilon}(x)$.
\begin{enumerate}
\item{
We assume that $V^{l,\epsilon},V^{u,\epsilon}\in \mathcal{C}^{3}(\mathbb{R})$ for
every fixed $\epsilon>0$ and that $V^{\epsilon}(x)\downarrow 0$ as
$\epsilon\downarrow 0$ (in particular $V^{0}(x)=0$). Moreover,
there exists a universal positive
constant $\zeta$ such that
\begin{equation}
\frac{1}{\epsilon}V^{\epsilon}(x)>\zeta>0 \textrm{ for every }
x\in\mathbb{R} \textrm{ and for every }
\epsilon>0. \label{DefinitionOfCrossSections1}
\end{equation} }
\item{We assume that the functions
\begin{eqnarray}
u(x)&:=&\lim_{\epsilon \downarrow 0 }u^{\epsilon}(x) \textrm{, } x\in\mathbb{R}\nonumber \\
v(x)&:=&\lim_{\epsilon \downarrow 0 }v^{\epsilon}(x) \textrm{, } x\in\mathbb{R}\setminus \lbrace 0\rbrace ,\label{uANDvForLimitingProcessIntro}
\end{eqnarray}
are well defined and  the limiting function $u(x)$ is continuous and strictly increasing whereas the limiting function $v(x)$ is right continuous and strictly increasing. In general, the function $u(x)$ can have countable  many
points where it is not differentiable and the function $v(x)$ can have countable
many points where it is not continuous or not differentiable. However,
here we assume for brevity that the only  non smoothness point is $x=0$. In other words, we assume that for $x\in\mathbb{R}\setminus \lbrace 0\rbrace$
\begin{equation}
V(x)=\frac{\partial V^{\epsilon}(x)}{\partial \epsilon}|_{\epsilon=0}>0,\label{DefinitionOfCrossSections2}
\end{equation}
and that the function $V(x)$ is smooth for $x\in\mathbb{R}\setminus \lbrace 0\rbrace$.

In addition, we assume that the first three derivatives of $V^{l,\epsilon}(x)$ and $V^{u,\epsilon}(x)$ (and in consequence of $V^{\epsilon}(x)$ as well) behave nicely for $|x|>0$ and for $\epsilon$ small. In particular, we assume that for any connected subset $K$ of $\mathbb{R}$ that is away from an arbitrarily small neighborhood of $x=0$ and for $\epsilon$ sufficiently small
\begin{equation}
|V^{l,\epsilon}_{x}(x)|+|V^{l,\epsilon}_{xx}(x)|+|V^{l,\epsilon}_{xxx}(x)|+|V^{u,\epsilon}_{x}(x)|+|V^{u,\epsilon}_{xx}(x)|+|V^{u,\epsilon}_{xxx}(x)|\leq C_{0}\epsilon \label{DefinitionOfCrossSections2a}
\end{equation}
uniformly in $x\in K$. Here, $C_{0}$ is a constant.

After
the proof of the main theorem (at the end of section 2), we mention the result for the case where
there exist more than one non smoothness point.
}
\item{Let $g^{\epsilon}(x)$ be a smooth function and let us define the quantity
\begin{displaymath}
\xi^{\epsilon}(g^{\epsilon}):=\sup_{|x|\leq 1}[|\frac{1}{\epsilon} [g^{\epsilon}_{x}(x)]^{3}|+|g^{\epsilon}_{x}(x)g^{\epsilon}_{xx}(x)|+|\epsilon g^{\epsilon}_{xxx}(x)|]
\end{displaymath}
We assume the following growth condition
\begin{equation}
\xi^{\epsilon}:=\xi^{\epsilon}(V^{\epsilon})+\xi^{\epsilon}(V^{l,\epsilon})+\xi^{\epsilon}(V^{u,\epsilon})\downarrow 0,  \textrm{ as } \epsilon\downarrow 0.\label{DefinitionOfCrossSections5}
\end{equation}
}
\end{enumerate}
\begin{rem}
Condition (\ref{DefinitionOfCrossSections5}), i.e.
$\xi^{\epsilon}\downarrow 0$, basically says that the behavior of $V^{l,\epsilon}(x)$ and $V^{u,\epsilon}(x)$ in the neighborhood of $x=0$ can be at most equally bad as described by $\xi^{\epsilon}(\cdot)$ for $\epsilon$ small. This condition will be used in the proof
of Lemma \ref{Lemma24} in section 4. Lemma 2.4 is essential for the proof of our main result. In particular, it provides us with the estimate of the expectation of the time it takes for the solution to (\ref{StochasticProcessWithReflection1}) to leave the neighborhood of the point $0$, as $\epsilon\downarrow 0$.
At the present moment, we do not know if this condition can be improved and this is subject to further research.
\end{rem}
In this paper we prove that under assumptions (\ref{DefinitionOfCrossSections1})-(\ref{DefinitionOfCrossSections5}), the $X^{\epsilon}_{t}$ component
of the process $(X^{\epsilon}_{t},Y^{\epsilon}_{t})$ converges weakly to a
one-dimensional strong Markov process, continuous with probability one. It behaves like a standard diffusion
process away from $0$ and has to satisfy a gluing condition at the point
of discontinuity $0$ as $\epsilon\downarrow 0$. More precisely,
we prove the following Theorem:

\vspace{0.2cm}

\begin{thm}
Let us assume that
(\ref{DefinitionOfCrossSections1})-(\ref{DefinitionOfCrossSections5})
hold. Let $X$ be the solution to the martingale problem for
\begin{equation}
A=\lbrace(f,Lf):f\in \mathcal{D} (A)\rbrace \label{LimitingProcess1}
\end{equation}
with
\begin{equation}
Lf(x)=D_{v}D_{u}f(x)\label{LimitingOperator0}
\end{equation}
and
\begin{eqnarray}
\mathcal{D} (A)=\lbrace &f:& f\in \mathcal{C}_{c}(\mathbb{R})\textrm{, with }f_{x},f_{xx} \in \mathcal{C}(\mathbb{R}\setminus\lbrace 0\rbrace),\nonumber\\
& & [u'(0+)]^{-1}f_{x}(0+)-[u'(0-)]^{-1} f_{x}(0-)=[v(0+)-v(0-)] Lf(0) \nonumber\\
& &  \textrm{ and } Lf(0)=\lim_{x\rightarrow
0^{+}}Lf(x)=\lim_{x\rightarrow 0^{-}}Lf(x) \rbrace. \label{LimitingOperator0Condition1}
\end{eqnarray}

Then we have
\begin{equation}
X^{\epsilon}_{\cdot}\longrightarrow X_{\cdot}\textrm{ weakly in } \mathcal{C}_{0T}, \textrm{ for any } T<\infty, \textrm{ as } \epsilon\downarrow 0,\label{Claim1}
\end{equation}
where $\mathcal{C}_{0T}$ is the space of continuous functions in $[0,T]$.
\begin{flushright}
$\square$
\end{flushright}\label{Theorem11}
\end{thm}
As proved in Feller \cite{Feller} the martingale problem for $A$,  (\ref{LimitingProcess1}), has a unique
solution $X$. It is an asymmetric
Markov process with delay at the point of discontinuity $0$. In
particular, the asymmetry is due to the possibility of having $u'(0+)\neq u'(0-)$ (see
Lemma \ref{Lemma25}) whereas the delay  is because of the possibility of having $v(0+)\neq v(0-)$ (see Lemma \ref{Lemma24}).

For the convenience of the reader, we briefly recall the Feller characterization
of all one-dimensional Markov processes, that are continuous with
probability one (for more details see \cite{Feller}; also \cite{M1}).
All one-dimensional strong Markov processes, that are continuous with
probability one, can be characterized (under some
minimal regularity conditions) by a generalized second order
differential operator $D_{v}D_{u}f$ with respect to two increasing
functions $u(x)$ and $v(x)$; $u(x)$ is continuous, $v(x)$ is right
continuous. In addition, $D_{u}$, $D_{v}$ are differentiation
operators with respect to $u(x)$ and $v(x)$ respectively, which are
defined as follows:

$D_{u}f(x)$ exists if $D_{u}f(x+)=D_{u}f(x-)$, where the left derivative of $f$ with respect to $u$ is defined as follows:
\begin{displaymath}
D_{u}f(x-)=\lim_{h\downarrow 0}\frac{f(x-h)-f(x)}{u(x-h)-u(x)} \hspace{0.2cm} \textrm{ provided the limit exists.}
\end{displaymath}
The right derivative $D_{u}f(x+)$ is defined similarly. If $v$ is discontinuous at $y$ then
\begin{displaymath}
 D_{v}f(y)=\lim_{h\downarrow 0}\frac{f(y+h)-f(y-h)}{v(y+h)-v(y-h)}.
\end{displaymath}
A more detailed description of these Markov processes can be found in \cite{Feller} and \cite{M1}.

\begin{rem}
Notice that if the limit of $\frac{1}{\epsilon}V^{\epsilon}(x)$, as $\epsilon\downarrow 0$, is a smooth function then the limiting process $X$ described by Theorem \ref{Theorem11} coincides with (\ref{LimitingStochasticProcess1}).
\end{rem}
We conclude the introduction with a useful example.
Let us assume that $V^{\epsilon}(x)$ can be decomposed in three terms
\begin{equation}
V^{\epsilon}(x)=V^{\epsilon}_{1}(x)+V^{\epsilon}_{2}(x)+V^{\epsilon}_{3}(x),\label{DefinitionOfCrossSections1special}
\end{equation}
where the functions $V^{\epsilon}_{i}(x)$, for $i=1,2,3$, satisfy the following conditions:
\begin{enumerate}
\item{There exists a strictly positive, smooth function $V_{1}(x)>0$ such that
\begin{equation}
\frac{1}{\epsilon}V^{\epsilon}_{1}(x)\rightarrow V_{1}(x), \textrm{ as } \epsilon\downarrow 0,\label{DefinitionOfCrossSections2special}
\end{equation}
uniformly in $x\in \mathbb{R}$.}
\item{There exists a nonnegative constant $\beta\geq 0$ such that
\begin{equation}
\frac{1}{\epsilon}V^{\epsilon}_{2}(x)\rightarrow\beta \chi_{\{x>0\}}, \textrm{ as } \epsilon\downarrow 0,\label{DefinitionOfCrossSections3special}
\end{equation}
uniformly for every connected subset of $\mathbb{R}$ that is away from an arbitrary small neighborhood of $0$ and weakly within a neighborhood of $0$. Here $\chi_{A}$ is the indicator function of the set $A$.}
\item{
\begin{equation}
\frac{1}{\epsilon}V^{\epsilon}_{3}(x)\rightarrow
\mu \delta_{0}(x),  \textrm{ as } \epsilon\downarrow
0,\label{DefinitionOfCrossSections4special}
\end{equation}
in the weak sense. Here $\mu$ is a nonnegative constant and $\delta_{0}(x)$ is the  Dirac delta
distribution at $0$.}
\item{Condition (\ref{DefinitionOfCrossSections5}) holds.
}
\end{enumerate}
Let us define $\alpha=V_{1}(0)$. In this case the operator (\ref{LimitingOperator0}) and its domain of definition (\ref{LimitingOperator0Condition1}) for the limiting process $X$ become
\begin{eqnarray}
Lf(x)=\cases{\frac{1}{2}f_{xx}(x)+\frac{1}{2}\frac{d}{dx}[\ln (V_{1}(x) )] f_{x}(x), & $x<0$ \cr
                                   \frac{1}{2}f_{xx}(x)+\frac{1}{2}\frac{d}{dx}[\ln (V_{1}(x)+\beta)] f_{x}(x),& $x>0$, \cr
                                   }\label{LimitingOperator0special}
\end{eqnarray}
and
\begin{eqnarray}
\mathcal{D} (A)=\lbrace &f:& f\in \mathcal{C}_{c}(\mathbb{R})\textrm{, with }f_{x},f_{xx} \in \mathcal{C}(\mathbb{R}\setminus\lbrace 0\rbrace), \nonumber\\
& &  [\alpha+\beta]f_{x}(0+)-\alpha f_{x}(0-)=[2\mu] Lf(0) \label{LimitingOperator0Condition1special}\\
& &  \textrm{ and } Lf(0)=\lim_{x\rightarrow
0^{+}}Lf(x)=\lim_{x\rightarrow 0^{-}}Lf(x) \rbrace. \nonumber
\end{eqnarray}
For instance, consider $0<\delta=\delta(\epsilon)\ll 1$ a small $\epsilon-$dependent positive number and assume that:
\begin{enumerate}
\item{$V^{l,\epsilon}(x)=0$ and $V^{u,\epsilon}(x)=V^{\epsilon}(x)=V^{\epsilon}_{1}(x)+V^{\epsilon}_{2}(x)+V^{\epsilon}_{3}(x)$ where}
\item{$V^{\epsilon}_{1}(x)=\epsilon V_{1}(x)$, where $V_{1}(x)$ is any smooth, strictly positive function,}
\item{$V^{\epsilon}_{2}(x)=\epsilon V_{2}(\frac{x}{\delta})$, such that $V_{2}(\frac{x}{\delta})\rightarrow\beta \chi_{\{x>0\}}$, as  $\epsilon\downarrow 0$}
\item{$V^{\epsilon}_{3}(x)=\frac{\epsilon}{\delta}V_{3}(\frac{x}{\delta})$, such that $\frac{1}{\delta}V_{3}(\frac{x}{\delta})\rightarrow\mu \delta_{0}(x)$, as  $\epsilon\downarrow 0$}
\item{and  with $\delta$ chosen such that $\frac{\epsilon}{\delta^{3}}\downarrow 0$ as $\epsilon\downarrow 0$.}
\end{enumerate}
Then, it can be easily verified that (\ref{DefinitionOfCrossSections2special})-(\ref{DefinitionOfCrossSections4special}) and (\ref{DefinitionOfCrossSections5}) are satisfied. Moreover, in this case, we have $\mu=\int_{-\infty}^{\infty}V_{3}(x)dx$.

In section 2 we prove our main result assuming that we have all
needed estimates. After the proof of Theorem \ref{Theorem11}, we
state the result in the case that $\lim_{\epsilon\downarrow
0}\frac{1}{\epsilon}V^{\epsilon}(x)$ has more than one point of
discontinuity (Theorem \ref{Theorem21}). In section 3 we prove
relative compactness of $X^{\epsilon}_{t}$ (this follows basically
from \cite{FW3}) and we consider what happens outside a small
neighborhood of $x=0$. In section 4 we estimate the expectation of
the time it takes for the solution to
(\ref{StochasticProcessWithReflection1}) to leave the neighborhood
of the point $0$. The derivation of this estimate uses assumption
(\ref{DefinitionOfCrossSections5}). In section 5 we: $(a)$ prove
that the behavior of the process after it reaches $x=0$ does not
depend on where it came from, and $(b)$ calculate the limiting exit
probabilities of $(X^{\epsilon}_{t},Y^{\epsilon}_{t})$, from the
left and from the right, of a small neighborhood of $x=0$. The
derivation of these estimates is composed of two main ingredients.
The first one is the characterization of all one-dimensional Markov
processes, that are continuous with probability one, by generalized
second order operators introduced by Feller (see \cite{Feller}; also
\cite{M1}). The second one is a result of Khasminskii on invariant
measures \cite{K1}.

Lastly, we would like to mention here that one can similarly
consider narrow tubes, i.e. $y\in D^{\epsilon}_{x}\subset
\mathbb{R}^{n}$ for $n>1$, and prove a result similar to Theorem
\ref{Theorem11}.

\section{Proof of the Main Theorem}
Before proving Theorem \ref{Theorem11} we introduce some notation
and formulate the necessary lemmas. The lemmas are proved
in sections 3 to 5.

In this and the following sections we will denote by $C_{0}$ any
unimportant constants that do not depend on any small parameter. The
constants may change from place to place though, but they will
always be denoted by the same $C_{0}$.

For any  $B\subset \mathbb{R}$, we define the Markov time $\tau(B)=\tau^{\epsilon}_{x,y}(B)$ to be:
\begin{equation}
\tau^{\epsilon}_{x,y}(B)= \inf \lbrace t>0: X_{t}^{\epsilon,x,y}
\notin B\rbrace.
\end{equation}
Moreover, for $\kappa>0$, the term
$\tau^{\epsilon}_{x,y}(\pm\kappa)$ will denote the Markov time
$\tau^{\epsilon}_{x,y}(-\kappa,\kappa)$. In addition,
$\mathbb{E}^{\epsilon}_{x,y}$ will denote the expected value
associated with the probability measure
$\mathbb{P}^{\epsilon}_{x,y}$ that is induced by the process
$(X^{\epsilon,x,y}_{t},Y^{\epsilon,x,y}_{t})$.

For the sake of notational convenience we define the operators
\begin{eqnarray}
L_{-}f(x)&=&D_{v}D_{u}f(x) \textrm{ for } x<0\nonumber\\
L_{+}f(x)&=&D_{v}D_{u}f(x) \textrm{ for } x>0\label{LimitingOperator1}
\end{eqnarray}
Furthermore, when we write $(x,y)\in A\times B$ we will mean
$(x,y)\in\lbrace (x,y): x\in A, y\in B\rbrace$.

Most of the processes, Markov times and sets that will be mentioned
below will depend on $\epsilon$. For notational convenience however,
we shall not always incorporate this dependence into the notation.
So the reader should be careful to distinguish between objects that
depend and do not depend on $\epsilon$.

Throughout this paper
$0<\kappa_{0}<\kappa$ will be small positive constants.
We may not always mention the relation
between these parameters but we will always assume it. Moreover $\kappa_{\eta}$ will denote a small positive number that depends on another small positive number $\eta$.

Moreover, one can write down the normal vector $\gamma^{\epsilon}(x,y)$ explicitly:
\begin{eqnarray}
\gamma^{\epsilon}(x,y)=\cases{\frac{1}{\sqrt{1+[V^{u,\epsilon}_{x}(x)]^{2}}}(V^{u,\epsilon}_{x}(x),-1), & $y=V^{u,\epsilon}(x)$ \cr
                                   \frac{1}{\sqrt{1+[V^{l,\epsilon}_{x}(x)]^{2}}}(V^{l,\epsilon}_{x}(x),1), & $y=-V^{l,\epsilon}(x)$. \cr
                                   }\nonumber
\end{eqnarray}


\vspace{0.2cm}

\begin{lem}
For any $(x,y) \in D^{\epsilon}$, let  $\mathbb{Q}^{\epsilon}_{x,y}$ be the
family of distributions of
$X^{\epsilon}_{\cdot}$ in the space $\mathcal{C}[0,\infty)$ of continuous functions $[0,\infty)\rightarrow \mathbb{R}$ that correspond to the probabilities $\mathbb{P}^{\epsilon}_{x,y}$. Assume that for any $(x,y) \in D^{\epsilon}$ the family of distributions $\mathbb{Q}^{\epsilon}_{x,y}$  for all $\epsilon\in(0,1)$
is tight. Moreover suppose that for any compact set $K\subseteq \mathbb{R}$, any function $f\in \mathcal{D}
(A)$ and for every $\lambda>0$ we have:
\begin{equation}
\mathbb{E}^{\epsilon}_{x,y} \int_{0}^{\infty}e^{-\lambda t} [\lambda
f(X^{\epsilon}_{t})-Lf(X^{\epsilon}_{t})]dt - f(x)\rightarrow
0,\label{MainCondition1}
\end{equation}
as $\epsilon \downarrow 0$, uniformly in $(x,y)\in K\times D^{\epsilon}_{x}$.

The measures $\mathbb{Q}^{\epsilon}_{x,y}$  corresponding to $\mathbb{P}^{\epsilon}_{x,y}$ converge weakly to the probability measure $\mathbb{P}_{x}$ that is induced by $X_{\cdot}$ as $\epsilon \downarrow 0$.
\begin{flushright}
$\square$
\end{flushright}\label{Lemma21}
\end{lem}

\begin{lem} The family of distributions
$\mathbb{Q}^{\epsilon}_{x,y}$ in the space of continuous functions $[0,\infty)\rightarrow \mathbb{R}$ corresponding to $\mathbb{P}^{\epsilon}_{x,y}$  for small nonzero $\epsilon$ is tight.
\begin{flushright}
$\square$
\end{flushright}\label{Lemma22}
\end{lem}
\begin{lem}
Let  $0<x_{1}<x_{2}$ be fixed positive numbers and $f$ be a three times continuously
differentiable function in $[ x_{1}, x_{2}]$. Then for every
$\lambda>0$:
\begin{equation}
\mathbb{E}^{\epsilon}_{x,y}[e^{-\lambda \tau^{\epsilon}( x_{1},
x_{2})}f(X^{\epsilon}_{\tau^{\epsilon}( x_{1},
x_{2})})+\int_{0}^{\tau^{\epsilon}( x_{1}, x_{2})}e^{-\lambda t}
[\lambda f(X^{\epsilon}_{t})-Lf(X^{\epsilon}_{t})]dt]\rightarrow
f(x),\nonumber
\end{equation}
as $\epsilon \downarrow 0$, uniformly in $(x,y)$ such that $(x,y)\in
[x_{1},x_{2}]\times  D^{\epsilon}_{x}$. The statement holds true for
$\tau^{\epsilon}(-x_{2},-x_{1})$ in place of $\tau^{\epsilon}(
x_{1}, x_{2})$ as well.
\begin{flushright}
$\square$
\end{flushright}\label{Lemma23}
\end{lem}
\begin{lem} Define $\theta=\frac{v(0+)-v(0-)}{[u'(0+)]^{-1}+[u'(0-)]^{-1}}$. For every $\eta>0$ there exists a $\kappa_{\eta}>0$ such that for every $0<\kappa<\kappa_{\eta}$ and for sufficiently small $\epsilon$
\begin{equation}
|\mathbb{E}^{\epsilon}_{x,y}\tau^{\epsilon}(\pm \kappa) -\kappa \theta|\leq\kappa \eta, \nonumber
\end{equation}
for all $(x,y)\in [-\kappa,  \kappa]\times  D^{\epsilon}_{x}$.
Here, $\tau^{\epsilon}(\pm \kappa)=\tau^{\epsilon}(-\kappa, \kappa)$ is the exit time from the
interval $(-\kappa, \kappa)$.
\begin{flushright}
$\square$
\end{flushright}\label{Lemma24}
\end{lem}

\begin{lem} Define $p_{+}=\frac{[u'(0+)]^{-1}}{[u'(0+)]^{-1}+[u'(0-)]^{-1}}$ and
$p_{-}=\frac{[u'(0-)]^{-1}}{[u'(0+)]^{-1}+[u'(0-)]^{-1}}$. For every
$\eta>0$ there exists a $\kappa_{\eta}>0$ such that for every $0<\kappa<\kappa_{\eta}$
there exists a positive $\kappa_{0}=\kappa_{0}(\kappa)$ such that for sufficiently small $\epsilon$
\begin{eqnarray}
|\mathbb{P}^{\epsilon}_{x,y}\lbrace
X^{\epsilon}_{\tau^{\epsilon}(\pm \kappa)}=\kappa \rbrace-p_{+}|&\leq&
\eta, \nonumber\\
|\mathbb{P}^{\epsilon}_{x,y}\lbrace
X^{\epsilon}_{\tau^{\epsilon}(\pm \kappa)}=-\kappa \rbrace-p_{-}|&\leq&
\eta, \nonumber
\end{eqnarray}
for all $(x,y)$ such that $(x,y)\in [-\kappa_{0},  \kappa_{0}]\times
D^{\epsilon}_{x}$.
\begin{flushright}
$\square$
\end{flushright}\label{Lemma25}
\end{lem}

\vspace{0.2cm}

\begin{proof}[Proof of Theorem \ref{Theorem11}]

We will make use of Lemma \ref{Lemma21}. The tightness required in
Lemma \ref{Lemma21} is the statement of Lemma \ref{Lemma22}. Thus it
remains to prove that (\ref{MainCondition1}) holds.

Let $\lambda>0$, $(x,y) \in D^{\epsilon}$ and $f\in \mathcal{D}(A)$
be fixed.  In addition let $\eta>0$ be
an arbitrary positive number.

Choose $0<x_{*}<\infty$ so that
\begin{equation}
\mathbb{E}^{\epsilon}_{x,y}e^{-\lambda \tilde{\tau}^{\epsilon}}<\frac{\eta}{\Vert f \Vert+\lambda^{-1}\Vert\lambda f-Lf\Vert} \label{DefinitionForT}
\end{equation}
for sufficiently small $\epsilon$, where $\tilde{\tau}^{\epsilon}=\inf\lbrace t>0:|X^{\epsilon}_{t}|\geq x_{*}\rbrace$.
It is Lemma \ref{Lemma22} that makes such a choice possible. We assume that $x_{*}>|x|$.

To prove (\ref{MainCondition1}) it is enough to show that for every
$\eta>0$ there exists an $\epsilon_{0}>0$, independent of $(x,y)$, such that for every
$0<\epsilon<\epsilon_{0}$:
\begin{equation}
|\mathbb{E}^{\epsilon}_{x,y}[e^{-\lambda \tilde{\tau}}f(X^{\epsilon}_{\tilde{\tau}})-f(x)+\int_{0}^{\tilde{\tau}}e^{-\lambda t} [\lambda f(X^{\epsilon}_{t})-Lf(X^{\epsilon}_{t})]dt] |<\eta.\label{GoalToProve}
\end{equation}
Choose $\epsilon$  small and  $0<\kappa_{0}<\kappa$ small positive numbers. We consider two cycles of
Markov times $\lbrace \sigma_{n}\rbrace$ and $\lbrace
\tau_{n}\rbrace$ such that:
\begin{displaymath}
0=\sigma_{0}\leq\tau_{1}\leq\sigma_{1}\leq\tau_{2}\leq\dots
\end{displaymath}
where:
\begin{eqnarray}
\tau_{n}&=&\tilde{\tau}\wedge\inf \lbrace t>\sigma_{n-1}: |X^{\epsilon}_{t}|\geq\kappa \rbrace \nonumber\\
\sigma_{n}&=&\tilde{\tau}\wedge\inf \lbrace t>\tau_{n}: |X^{\epsilon}_{t}|\in \{ \kappa_{0}, x_{*} \} \rbrace \label{MarkovTimes1}
\end{eqnarray}
In figure 1 we see a trajectory of the process $Z_{t}^{\epsilon}=(X_{t}^{\epsilon},Y_{t}^{\epsilon})$ along with its associated Markov chains $\{Z_{\sigma_{n}}^{\epsilon}\}$ and $\{Z_{\tau_{n}}^{\epsilon}\}$. We will also write $z=(x,y)$ for the initial point.

\begin{figure}[ht]
\begin{center}
\includegraphics[scale=1, width=9 cm, height=7 cm]{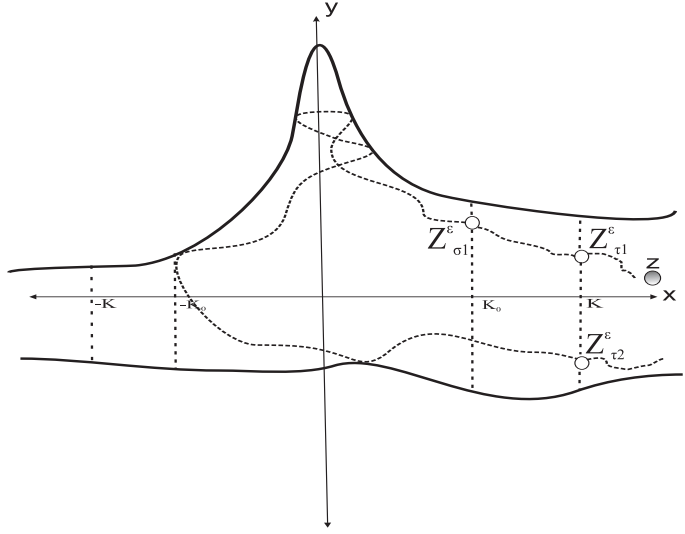}
\caption{ $z=(x,y)$ is the initial point and
$Z_{t}^{\epsilon}=(X_{t}^{\epsilon},Y_{t}^{\epsilon})$.}
\end{center}
\end{figure}

We denote by $\chi_{A}$ the indicator function of the set $A$. The
difference in (\ref{GoalToProve}) can be represented as the sum over
time intervals from $\sigma_{n}$ to $\tau_{n+1}$ and from $\tau_{n}$
to $\sigma_{n}$. It is equal to:

\begin{eqnarray}
& &\mathbb{E}^{\epsilon}_{z}[e^{-\lambda \tilde{\tau}}f(X^{\epsilon}_{\tilde{\tau}})-f(x)+\int_{0}^{\tilde{\tau}}e^{-\lambda t} [\lambda f(X^{\epsilon}_{t})-Lf(X^{\epsilon}_{t})]dt]=\label{GoalToProve02}\\
&=&\mathbb{E}^{\epsilon}_{z}\sum_{n=0}^{\infty}[e^{-\lambda
\tau_{n+1}^{\epsilon}}f(X^{\epsilon}_{\tau_{n+1}^{\epsilon}})-e^{-\lambda
\sigma_{n}^{\epsilon}}f(X_{\sigma_{n}^{\epsilon}})
+\int_{\sigma_{n}^{\epsilon}}^{\tau_{n+1}^{\epsilon}}e^{-\lambda t} [\lambda f(X^{\epsilon}_{t})-Lf(X^{\epsilon}_{t})]dt]]+\nonumber\\
&+&\mathbb{E}^{\epsilon}_{z}\sum_{n=1}^{\infty}[e^{-\lambda
\sigma_{n}^{\epsilon}}f(X^{\epsilon}_{\sigma_{n}^{\epsilon}})-e^{-\lambda
\tau_{n}^{\epsilon}}f(X_{\tau_{n}^{\epsilon}})+\int_{\tau_{n}^{\epsilon}}^{\sigma_{n}^{\epsilon}}e^{-\lambda
t} [\lambda f(X^{\epsilon}_{t})-Lf(X^{\epsilon}_{t})]dt]\nonumber
\end{eqnarray}
The formally infinite sums are finite for every trajectory for which
$\tilde{\tau}<\infty$. Assuming that we can write the expectation of
the infinite sums as the infinite sum of the expectations, the
latter equality becomes
\begin{eqnarray}
& &\mathbb{E}^{\epsilon}_{z}[e^{-\lambda \tilde{\tau}}f(X^{\epsilon}_{\tilde{\tau}})-f(x)+\int_{0}^{\tilde{\tau}}e^{-\lambda t} [\lambda f(X^{\epsilon}_{t})-Lf(X^{\epsilon}_{t})]dt]=\label{GoalToProve2}\\
&=&\sum_{n=0}^{\infty}\mathbb{E}^{\epsilon}_{z}[e^{-\lambda
\tau_{n+1}^{\epsilon}}f(X^{\epsilon}_{\tau_{n+1}^{\epsilon}})-e^{-\lambda
\sigma_{n}^{\epsilon}}f(X_{\sigma_{n}^{\epsilon}})
+\int_{\sigma_{n}^{\epsilon}}^{\tau_{n+1}^{\epsilon}}e^{-\lambda t} [\lambda f(X^{\epsilon}_{t})-Lf(X^{\epsilon}_{t})]dt]]+\nonumber\\
&+&\sum_{n=1}^{\infty}\mathbb{E}^{\epsilon}_{z}[e^{-\lambda
\sigma_{n}^{\epsilon}}f(X^{\epsilon}_{\sigma_{n}^{\epsilon}})-e^{-\lambda
\tau_{n}^{\epsilon}}f(X_{\tau_{n}^{\epsilon}})+\int_{\tau_{n}^{\epsilon}}^{\sigma_{n}^{\epsilon}}e^{-\lambda
t} [\lambda f(X^{\epsilon}_{t})-Lf(X^{\epsilon}_{t})]dt]\nonumber
\end{eqnarray}
The aforementioned calculation can be done if
$$\sum_{n=0}^{\infty}\mathbb{E}^{\epsilon}_{z}[e^{-\lambda
\sigma_{n}}\chi_{\sigma_{n}<\tilde{\tau}}],\sum_{n=1}^{\infty}\mathbb{E}^{\epsilon}_{z}[e^{-\lambda
\tau_{n}}\chi_{\tau_{n}<\tilde{\tau}}]<\infty.$$ Indeed, by Markov
property we have:
\begin{eqnarray}
\mathbb{E}^{\epsilon}_{z}[e^{-\lambda
\sigma_{n}}\chi_{\sigma_{n}<\tilde{\tau}}]&\leq&\mathbb{E}^{\epsilon}_{z}[e^{-\lambda
\tau_{n}}\chi_{\tau_{n}<\tilde{\tau}}]\max_{|x|=\kappa, y\in
D^{\epsilon}_{x}} \phi_{1}^{\epsilon}(x,y)\leq \nonumber\\
 &\leq&
\mathbb{E}^{\epsilon}_{z}[e^{-\lambda
\sigma_{n-1}}\chi_{\sigma_{n-1}<\tilde{\tau}}]\max_{|x|=\kappa, y\in
D^{\epsilon}_{x}} \phi_{1}^{\epsilon}(x,y),\label{InfiniteSums0}
\end{eqnarray}
where $\phi_{1}^{\epsilon}(x,y)=\mathbb{E}_{x,y}^{\epsilon}[e^{-\lambda \sigma_{1}}\chi_{|X^{\epsilon}_{\sigma_{1}}|=\kappa_{0}}]$. So by induction we have
\begin{equation}
\sum_{n=1}^{\infty}\mathbb{E}^{\epsilon}_{z}[e^{-\lambda
\tau_{n}}\chi_{\tau_{n}<\tilde{\tau}}]\leq
\sum_{n=0}^{\infty}\mathbb{E}^{\epsilon}_{z}[e^{-\lambda
\sigma_{n}}\chi_{\sigma_{n}<\tilde{\tau}}]\leq
\frac{1}{1-\max_{|x|=\kappa, y\in D^{\epsilon}_{x}}
\phi_{1}^{\epsilon}(x,y)}.\label{InfiniteSums1}
\end{equation}
Clearly  $\max_{|x|=\kappa, y\in D^{\epsilon}_{x}} \phi_{1}^{\epsilon}(x,y)<1$
for $\kappa \in (\kappa_{0},x_{*})$. Therefore equality (\ref{GoalToProve2}) is valid.

However we need to know how the sums in (\ref{InfiniteSums1}) behave in terms of $\kappa$.
To this end we apply Lemma 2.3 to the function $g$ that is the solution to
\begin{eqnarray}
\lambda g-L_{\pm} g &=&0 \hspace{0.5cm}\mathrm{ in }\hspace{0.2cm} x \in (\pm \kappa_{0},\pm x_{*})\nonumber\\
g(\pm \kappa_{0})&=& 1\nonumber\\
g(\pm x_{*})&=& 0 \label{MaximumPrinciple2}
\end{eqnarray}
By Lemma \ref{Lemma23} we know that $g(x)$ approximates $\phi_{1}^{\epsilon}(x,y)$ for $|x|\in [\kappa_{0},x_{*}]$ as $\epsilon \downarrow 0$. The idea is to bound $\phi_{1}^{\epsilon}(x,y)$ using $g(x)$. It follows by (\ref{MaximumPrinciple2}) (for more details see the related discussion in section $8.3$ of \cite{FW1}, page 306) that there exists a positive constant $C_{0}$ that is independent of $\epsilon$ and a positive constant $\kappa^{'}$ such that for every $\kappa<\kappa^{'}$ and for all $\kappa_{0}<\kappa_{0}^{'}(\kappa)$ we have
$g(\pm \kappa)\leq 1- C_{0} \kappa$.

So we conclude for $\epsilon$ and $\kappa$ sufficiently small that
\begin{equation}
\sum_{n=1}^{\infty}\mathbb{E}^{\epsilon}_{z}[e^{-\lambda \tau_{n}}\chi_{\tau_{n}<\tilde{\tau}}]\leq \sum_{n=0}^{\infty}\mathbb{E}^{\epsilon}_{z}[e^{-\lambda \sigma_{n}}\chi_{\sigma_{n}<\tilde{\tau}}]\leq\frac{C_{0}}{\kappa}.\label{InfiniteSums2}
\end{equation}
By the strong Markov property with respect to the Markov times
$\tau_{n}$ and $\sigma_{n}$ equality (\ref{GoalToProve2}) becomes
\begin{eqnarray}
& &\mathbb{E}^{\epsilon}_{z}[e^{-\lambda \tilde{\tau}}f(X^{\epsilon}_{\tilde{\tau}})-f(x)+\int_{0}^{\tilde{\tau}}e^{-\lambda t} [\lambda f(X^{\epsilon}_{t})-Lf(X^{\epsilon}_{t})]dt]=\nonumber\\
&=&\sum_{n=0}^{\infty}\mathbb{E}^{\epsilon}_{z}[e^{-\lambda
\sigma_{n}}\chi_{\sigma_{n}<\tilde{\tau}}
\mathbb{E}^{\epsilon}_{Z^{\epsilon}_{\sigma_{n}^{\epsilon}}}[e^{-\lambda \tau_{1}^{\epsilon}}f(X^{\epsilon}_{\tau_{1}^{\epsilon}})-f(X_{\sigma_{n}^{\epsilon}})+\int_{0}^{\tau_{1}^{\epsilon}}e^{-\lambda t} [\lambda f(X^{\epsilon}_{t})-Lf(X^{\epsilon}_{t})]dt]]+\nonumber\\
&+&\sum_{n=1}^{\infty}\mathbb{E}^{\epsilon}_{z}[e^{-\lambda
\tau_{n}}\chi_{\tau_{n}<\tilde{\tau}}\mathbb{E}^{\epsilon}_{Z^{\epsilon}_{\tau_{n}}}[e^{-\lambda
\sigma_{1}^{\epsilon}}f(X^{\epsilon}_{\sigma_{1}^{\epsilon}})-f(X_{\tau_{n}^{\epsilon}})+\int_{0}^{\sigma_{1}^{\epsilon}}e^{-\lambda
t} [\lambda f(X^{\epsilon}_{t})-Lf(X^{\epsilon}_{t})]dt]]\nonumber\\
&=&\sum_{n=0}^{\infty}\mathbb{E}^{\epsilon}_{z}[e^{-\lambda
\sigma_{n}}\chi_{\sigma_{n}<\tilde{\tau}}
\phi_{2}^{\epsilon}(Z^{\epsilon}_{\sigma_{n}^{\epsilon}})]+
\sum_{n=1}^{\infty}\mathbb{E}^{\epsilon}_{z}[e^{-\lambda
\tau_{n}}\chi_{\tau_{n}<\tilde{\tau}}\phi_{3}^{\epsilon}(Z^{\epsilon}_{\tau_{n}})]\label{GoalToProve2a}
\end{eqnarray}
where
\begin{equation}
\phi_{2}^{\epsilon}(x,y)=\mathbb{E}^{\epsilon}_{x,y}[e^{-\lambda
\tau_{1}}f(X^{\epsilon}_{\tau_{1}})-f(x)+\int_{0}^{\tau_{1}}e^{-\lambda
t} [\lambda f(X^{\epsilon}_{t})-Lf(X^{\epsilon}_{t})]dt]
\end{equation}
and
\begin{equation}
\phi_{3}^{\epsilon}(x,y)=\mathbb{E}^{\epsilon}_{x,y}[e^{-\lambda
\sigma_{1}}f(X^{\epsilon}_{\sigma_{1}})-f(x)+\int_{0}^{\sigma_{1}}e^{-\lambda
t} [\lambda f(X^{\epsilon}_{t})-Lf(X^{\epsilon}_{t})]dt]
\end{equation}

Because of (\ref{InfiniteSums2}), equality (\ref{GoalToProve2a})
becomes
\begin{eqnarray}
& &|\mathbb{E}^{\epsilon}_{x,y}[e^{-\lambda \tilde{\tau}}f(X^{\epsilon}_{\tilde{\tau}})-f(x)+\int_{0}^{\tilde{\tau}}e^{-\lambda t} [\lambda f(X^{\epsilon}_{t})-Lf(X^{\epsilon}_{t})]dt]|\leq\nonumber\\
&\leq&|\phi_{2}^{\epsilon}(x,y)| + \frac{C_{0}}{\kappa}\left[\max_{|x|=\kappa_{0}, y\in D^{\epsilon}_{x}} |\phi_{2}^{\epsilon}(x,y)|+\max_{|x|=\kappa, y\in D^{\epsilon}_{x}} |\phi_{3}^{\epsilon}(x,y)|\right] \label{GoalToProve3}
\end{eqnarray}

By Lemma \ref{Lemma23} we get that  $\max_{|x|=\kappa, y\in
D^{\epsilon}_{x}} |\phi_{3}^{\epsilon}(x,y)|$ is arbitrarily small
for sufficiently small $\epsilon$, so
\begin{equation}
\frac{C_{0}}{\kappa}\max_{|x|=\kappa, y\in D^{\epsilon}_{x}} |\phi_{3}^{\epsilon}(x,y)|\leq \frac{\eta}{3} \label{GoalToProve3Part1}
\end{equation}
Therefore, it remains to consider the terms
$|\phi_{2}^{\epsilon}(x,y)|$, where $(x,y)$ is the initial point, and $\frac{1}{\kappa}\max_{|x|=\kappa_{0}, y\in D^{\epsilon}_{x}}
|\phi_{2}^{\epsilon}(x,y)|$.

Firstly, we consider the term $|\phi_{2}^{\epsilon}(x,y)|$, where $(x,y)$ is the initial point. Clearly, if $|x|>\kappa$, then Lemma \ref{Lemma23} implies that  $|\phi_{2}^{\epsilon}(x,y)|$ is arbitrarily small
for sufficiently small $\epsilon$, so
\begin{equation}
|\phi_{2}^{\epsilon}(x,y)|\leq \frac{\eta}{3}. \label{GoalToProve3Part1a}
\end{equation}

We consider now the case $|x|\leq\kappa$. Clearly, in this case Lemma \ref{Lemma23} does not apply. However, one can use the continuity of $f$ and Lemma \ref{Lemma24}, as the following calculations show. We have:
\begin{displaymath}
|\phi_{2}^{\epsilon}(x,y)|\leq \mathbb{E}^{\epsilon}_{x,y}|f(X^{\epsilon}_{\tau_{1}})-f(x)|+|\lambda\Vert f \Vert+\Vert\lambda f-Lf\Vert|
\mathbb{E}^{\epsilon}_{x,y}\int_{0}^{\tau_{1}}e^{-\lambda t} dt.
\end{displaymath}
Choose now a positive $\kappa'$ so that
\begin{displaymath}
|f(x)-f(0)|<\frac{\eta}{6}\textrm{, for all } |x|\leq\kappa'
\end{displaymath}
and that
\begin{displaymath}
\mathbb{E}^{\epsilon}_{x,y}\int_{0}^{\tau_{1}}e^{-\lambda t} dt\leq \frac{\eta}{6[\lambda\Vert f \Vert+\Vert\lambda f-Lf\Vert]}
\end{displaymath}
for sufficiently small $\epsilon$ and for all $|x|\leq\kappa'$. Therefore, for $\kappa\leq\kappa'$ and for sufficiently small $\epsilon$ we have
\begin{equation}
|\phi_{2}^{\epsilon}(x,y)|\leq \frac{\eta}{3}\textrm{, for all }(x,y)\in [-\kappa,\kappa]\times D^{\epsilon}_{x}. \label{GoalToProve3Part1b}
\end{equation}

Secondly, we consider the term  $\frac{1}{\kappa}\max_{|x|=\kappa_{0}, y\in D^{\epsilon}_{x}}
|\phi_{2}^{\epsilon}(x,y)|$. Here, we need a sharper estimate because of the factor $\frac{1}{\kappa}$. We will prove that for $(x,y)\in \lbrace \pm
\kappa_{0} \rbrace \times D^{\epsilon}_{\pm \kappa_{0}}$ and for $\epsilon$ sufficiently small
\begin{equation}
|\phi_{2}^{\epsilon}(x,y)|\leq \kappa \frac{\eta}{3C_{0}}.\label{GoalToProve3Part2Step1a}
\end{equation}
For $(x,y)\in \lbrace \pm \kappa_{0} \rbrace \times
D^{\epsilon}_{\pm \kappa_{0}}$ we have
\begin{eqnarray}
|\phi_{2}^{\epsilon}(x,y)|&=& |\mathbb{E}^{\epsilon}_{x,y}[e^{-\lambda \tau_{1}}f(X^{\epsilon}_{\tau_{1}})-f(x)+\int_{0}^{\tau_{1}}e^{-\lambda t} [\lambda f(X^{\epsilon}_{t})-Lf(X^{\epsilon}_{t})]dt]|\nonumber\\
&\leq & |\mathbb{E}^{\epsilon}_{x,y}[f(X^{\epsilon}_{\tau_{1}})-f(x)-\kappa\theta Lf(0)]|+\nonumber\\
& & +  |\mathbb{E}^{\epsilon}_{x,y}[\tau_{1}Lf(0)- \int_{0}^{\tau_{1}}e^{-\lambda t} Lf(X^{\epsilon}_{t})dt]|+\nonumber\\
& & +  |\mathbb{E}^{\epsilon}_{x,y}[\kappa\theta Lf(0)- \tau_{1}Lf(0)|+\nonumber\\
& & +  |\mathbb{E}^{\epsilon}_{x,y}[e^{-\lambda
\tau_{1}}f(X^{\epsilon}_{\tau_{1}})-f(X^{\epsilon}_{\tau_{1}})+
\int_{0}^{\tau_{1}}e^{-\lambda t} \lambda f(X^{\epsilon}_{t})dt]|,
\label{GoalToProve3Part2Step1}
\end{eqnarray}
where $\theta=\frac{v(0+)-v(0-)}{[u'(0+)]^{-1}+[u'(0-)]^{-1}}$.

Since the one-sided derivatives of $f$ exist, we may choose, a positive $\kappa_{\eta}$ such that for every $0<\kappa_{1}\leq\kappa_{\eta}$
\begin{equation}
|\frac{f(w)-f(0)}{w}-f_{x}(0+)|,|\frac{f(-w)-f(0)}{-w}-f_{x}(0-)|\leq \frac{\eta}{C_{0}},\label{GoalToProve3Part2Step2}
\end{equation}
for all $w\in (0,\kappa_{1})$.

Furthermore, by Lemma \ref{Lemma25} we can choose for sufficiently small
$\kappa_{2}>0$, a $\kappa_{0}(\kappa_{2})\in(0,\kappa_{2})$ such that for sufficiently small $\epsilon$
\begin{equation}
|\mathbb{P}^{\epsilon}_{x,y}\lbrace X^{\epsilon}_{\tau_{1}^{\epsilon}(\pm \kappa_{2})}=\pm \kappa_{2} \rbrace-p_{\pm}|\leq \frac{\eta}{C_{0}}  \label{GoalToProve3Part2Step3}
\end{equation}
for all $(x,y)$ such that $(x,y)\in [-\kappa_{0},\kappa_{0}]\times D^{\epsilon}_{x}$.

In addition, by Lemma \ref{Lemma24} we can choose for sufficiently small
$\kappa_{\eta}>0$,  a $\kappa_{3}\in(0,\kappa_{\eta})$ such that for sufficiently small $\epsilon$
\begin{equation}
|\mathbb{E}^{\epsilon}_{x,y}\tau_{1}^{\epsilon}(\pm\kappa_{3})-\kappa_{3}\theta|\leq \kappa_{3}\frac{\eta}{C_{0}}  \label{GoalToProve3Part2Step3a}
\end{equation}
for all  $(x,y)\in [-\kappa_{3},\kappa_{3}]\times D^{\epsilon}_{x}$.

Choose now $0<\kappa \leq \min\lbrace\kappa_{1},\kappa_{2},\kappa_{3}\rbrace$ and
$0<\kappa_{0} < \min\lbrace\kappa_{0}(\kappa_{2}),\kappa \rbrace$.

For sufficiently small
$\epsilon$ and for all $(x,y)\in \lbrace \pm \kappa_{0} \rbrace
\times  D^{\epsilon}_{\pm\kappa_{0}}$ we have
\begin{eqnarray}
& &
|\mathbb{E}^{\epsilon}_{x,y}f(X^{\epsilon}_{\tau_{1}})-f(x)-\kappa\theta
Lf(0)|\leq\nonumber\\
&\leq&|p_{+}[f(\kappa)-f(0)]+p_{-}[f(-\kappa)-f(0)]-\kappa \theta Lf(0)|+\nonumber\\
&+&|\mathbb{P}^{\epsilon}_{x,y}\lbrace
X^{\epsilon}_{\tau_{1}^{\epsilon}(\pm\kappa)}=\kappa
\rbrace-p_{+}||f(\kappa)-f(0)|+\nonumber\\
&+&|\mathbb{P}^{\epsilon}_{x,y}\lbrace
X^{\epsilon}_{\tau_{1}^{\epsilon}(\pm\kappa)}=-\kappa
\rbrace-p_{-}||f(-\kappa)-f(0)|+\nonumber\\
&+&|f(0)-f(x)|\label{GoalToProve3Part2Step4}
\end{eqnarray}
Because of (\ref{GoalToProve3Part2Step2}) and the gluing condition
$p_{+}f_{x}(0+)-p_{-}f_{x}(0-)=\theta Lf(0)$, the first summand on the right hand
side of (\ref{GoalToProve3Part2Step4}) satisfies
\begin{eqnarray}
& &|p_{+}[f(\kappa)-f(0)]+p_{-}[f(-\kappa)-f(0)]-\kappa \theta Lf(0)|\leq \nonumber\\
&\leq& |p_{+}\kappa f_{x}(0+)-p_{-}\kappa f_{x}(0-)-\kappa \theta Lf(0)|+\kappa \frac{\eta}{C_{0}}=\nonumber\\
&=& \kappa \frac{\eta}{C_{0}}.\label{GoalToProve3Part2Step5}
\end{eqnarray}
Moreover, for small enough $x\in \lbrace\pm \kappa_{0},\pm \kappa\rbrace$ we also have that
\begin{equation}
|f(x)-f(0)|\leq |x| |f_{x}(0 \pm)|+|x| \eta.\label{GoalToProve3Part2Step5a}
\end{equation}
The latter together with (\ref{GoalToProve3Part2Step3}) imply that for sufficiently small $\epsilon$
the second summand on the right hand side of
(\ref{GoalToProve3Part2Step4}) satisfies
\begin{equation}
|P^{\epsilon}_{x,y}\lbrace X^{\epsilon}_{\tau_{1}^{\epsilon}(\pm
\kappa)}=\kappa \rbrace-p_{+}||f(\kappa)-f(0)|\leq\kappa \frac{\eta}{C_{0}}.\label{GoalToProve3Part2Step6}
\end{equation}
A similar expression holds for the third summand on the right hand
side of (\ref{GoalToProve3Part2Step4}) as well. Therefore
(\ref{GoalToProve3Part2Step5})-(\ref{GoalToProve3Part2Step6}) and
the fact that we take $\kappa_{0}$ to be much smaller than $\kappa$ imply that for
all $(x,y)\in \lbrace \pm \kappa_{0} \rbrace \times
D^{\epsilon}_{\pm\kappa_{0}}$ and for $\epsilon$ sufficiently small, we have
\begin{equation}
|\mathbb{E}^{\epsilon}_{x,y}f(X^{\epsilon}_{\tau_{1}})-f(x)-\kappa \theta Lf(0)|\leq\kappa \frac{\eta}{C_{0}}.\label{GoalToProve3Part2Step7}
\end{equation}
The second term of the right hand side of
(\ref{GoalToProve3Part2Step1}) can also be bounded by $\kappa \frac{\eta}{C_{0}}$ for $\kappa$ and $\epsilon$ sufficiently small, as the following calculations show. For $(x,y)\in
\lbrace \pm \kappa_{0} \rbrace \times  D^{\epsilon}_{\pm\kappa_{0}}$ we have
\begin{eqnarray}
& & |\mathbb{E}^{\epsilon}_{x,y}[\tau_{1}Lf(0)- \int_{0}^{\tau_{1}}e^{-\lambda t} Lf(X^{\epsilon}_{t})dt]|\leq\nonumber\\
&\leq&|Lf(0)||\mathbb{E}^{\epsilon}_{x,y}[\tau_{1}-
\int_{0}^{\tau_{1}}e^{-\lambda t}dt]|+\sup_{|x|\leq
\kappa}|Lf(x)-Lf(0)|\mathbb{E}^{\epsilon}_{x,y}\tau_{1}\leq\label{GoalToProve3Part2Step8a}\\
&\leq&\lambda|Lf(0)|\mathbb{E}^{\epsilon}_{x,y}\tau_{1}\left[\sup_{(x,y)\in
\lbrace \pm \kappa_{0} \rbrace \times
D^{\epsilon}_{\pm\kappa_{0}}}\mathbb{E}^{\epsilon}_{x,y}\tau_{1}\right]+\sup_{|x|\leq
\kappa}|Lf(x)-Lf(0)|\mathbb{E}^{\epsilon}_{x,y}\tau_{1}\nonumber
\end{eqnarray}
Therefore, Lemma \ref{Lemma24} (in particular (\ref{GoalToProve3Part2Step3a})) and the continuity of the function
$Lf$ give us for $\kappa$ and $\epsilon$ sufficiently small that
\begin{equation}
|\mathbb{E}^{\epsilon}_{x,y}[\tau_{1}Lf(0)-
\int_{0}^{\tau_{1}}e^{-\lambda t} Lf(X^{\epsilon}_{t})dt]|\leq
\kappa \frac{\eta}{C_{0}}.\label{GoalToProve3Part2Step8}
\end{equation}
The third term of the right hand side of
(\ref{GoalToProve3Part2Step1}) is clearly bounded by $\kappa \frac{\eta}{C_{0}}$ for $\epsilon$ sufficiently small
by Lemma \ref{Lemma24}. As far as the fourth
term of the right hand side of (\ref{GoalToProve3Part2Step1}) is
concerned, one can use the continuity of $f$ together with
Lemma \ref{Lemma24}.

The latter, (\ref{GoalToProve3Part2Step7}), (\ref{GoalToProve3Part2Step8}) and (\ref{GoalToProve3Part2Step1}) finally
give us that
\begin{equation}
\frac{1}{\kappa}\max_{|x|=\kappa_{0}, y\in D^{\epsilon}_{x}} |\phi_{2}^{\epsilon}(x,y)]|\leq \frac{\eta}{3C_{0}}. \label{GoalToProve3Part2}
\end{equation}
Of course, the constants $C_{0}$ that appear in the relations above are not the same, but for notational convenience they are all denoted by the same symbol $C_{0}$.

So, we finally get by (\ref{GoalToProve3Part2}), (\ref{GoalToProve3Part1}), (\ref{GoalToProve3Part1a}), (\ref{GoalToProve3Part1b}) and (\ref{GoalToProve3}) that
\begin{equation}
|\mathbb{E}^{\epsilon}_{x,y}[e^{-\lambda \tilde{\tau}}f(X^{\epsilon}_{\tilde{\tau}})-f(x)+\int_{0}^{\tilde{\tau}}e^{-\lambda t} [\lambda f(X^{\epsilon}_{t})-Lf(X^{\epsilon}_{t})]dt]|\leq \eta.
\end{equation}
This concludes the proof of Theorem \ref{Theorem11}.
\end{proof}
In case $\lim_{\epsilon\downarrow
0}\frac{1}{\epsilon}V^{\epsilon}(x)$ has more than one points of
discontinuity, one can similarly prove the following theorem. Hence,
the limiting Markov process $X$ may be asymmetric at some point
$x_{1}$, have delay at some other point $x_{2}$ or have both
irregularities at another point $x_{3}$.
\begin{thm}
Let us assume that $\frac{1}{\epsilon}V^{\epsilon}(x)$ has a finite
number of discontinuities, as described by
(\ref{DefinitionOfCrossSections1})-(\ref{DefinitionOfCrossSections5}),
at $x_{i}$ for $i\in\lbrace 1,\cdots,m\rbrace$. Let $X$ be the
solution to the martingale problem for
\begin{displaymath}
A=\lbrace(f,Lf):f\in \mathcal{D} (A)\rbrace
\end{displaymath}
with
\begin{displaymath}
Lf(x)=D_{v}D_{u}f(x)\nonumber
\end{displaymath}
and
\begin{eqnarray}
\mathcal{D} (A)=\lbrace &f:& f\in \mathcal{C}_{c}(\mathbb{R})\textrm{, with }f_{x},f_{xx} \in \mathcal{C}(\mathbb{R}\setminus\lbrace x_{1},\cdots, x_{m}  \rbrace) \nonumber\\
& & [u'(x_{i}+)]^{-1}f_{x}(x_{i}+)-[u'(x_{i}-)]^{-1}f_{x}(x_{i}-)=[v(x_{i}+)-v(x_{i}-)] Lf(x_{i}) \nonumber\\
& &  \textrm{ and } Lf(x_{i})=\lim_{x\rightarrow x_{i}^{+}}Lf(x)=\lim_{x\rightarrow x_{i}^{-}}Lf(x)  \textrm{ for } i\in\lbrace1,\cdots,m\rbrace \rbrace, \nonumber
\end{eqnarray}
Then we have
\begin{displaymath}
X^{\epsilon}_{\cdot}\longrightarrow X_{\cdot}\textrm{ weakly in } \mathcal{C}_{0T}, \textrm{ for any } T<\infty, \textrm{ as } \epsilon\downarrow 0.
\end{displaymath}
\begin{flushright}
$\square$
\end{flushright}\label{Theorem21}
\end{thm}

\section{Proof of Lemmata \ref{Lemma21}, \ref{Lemma22} and \ref{Lemma23}}
\begin{proof}[Proof of Lemma \ref{Lemma21}]

The proof is very similar to the proof of  Lemma 8.3.1 in
\cite{FW1}, so it will not be repeated here.
\end{proof}
\begin{proof}[Proof of Lemma \ref{Lemma22}]

The tool that is used to establish tightness of
$\mathbb{P}^{\epsilon}$ is the martingale-problem approach of
Stroock-Varadhan \cite{SV1}. In particular we can apply Theorem 2.1
of \cite{FW3}. The proof is almost identical to the part of the
proof of  Theorem 6.1 in \cite{FW3} where pre-compactness is proven
for the Wiener process with reflection in narrow-branching tubes.
\end{proof}
Before proving Lemma \ref{Lemma23} we introduce the following
diffusion process. Let $\hat{X}_{t}^{\epsilon}$ be the
one-dimensional process that is the solution to:
\begin{equation}
\hat{X}_{t}^{\epsilon}= x+
W_{t}^{1}+\int_{0}^{t}\frac{1}{2}\frac{V^{\epsilon}_{x}(\hat{X}_{s}^{\epsilon})}{V^{\epsilon}(\hat{X}_{s}^{\epsilon})}ds,
\label{EscordProcess}
\end{equation}
where $V_{x}^{\epsilon}(x)=\frac{dV^{\epsilon}(x)}{dx}$. The
process $\hat{X}^{\epsilon}$ is solution to the martingale
problem for $\hat{A}^{\epsilon}=\lbrace(f,\hat{L}^{\epsilon}f):f\in
\mathcal{D} (\hat{A}^{\epsilon})\rbrace$ with
\begin{equation}
\hat{L}^{\epsilon}=\frac{1}{2}\frac{d^{2}}{dx^{2}}+\frac{1}{2}\frac{V^{\epsilon}_{x}(\cdot)}{V^{\epsilon}(\cdot)}\frac{d}{dx}. \label{OperatorOfHatProcess}
\end{equation}
and
\begin{equation}
\mathcal{D} (\hat{A}^{\epsilon})=\lbrace f: f\in \mathcal{C}^{2}_{c}(\mathbb{R})\rbrace \label{EscortProcessOperator}
\end{equation}
A simple calculation shows that
\begin{displaymath}
\hat{L}^{\epsilon}f(x)=D_{v^{\epsilon}}D_{u^{\epsilon}}f(x).
\end{displaymath}
where the $u^{\epsilon}(x)$ and $v^{\epsilon}(x)$ functions are defined  by (\ref{uANDvForEscortProcessIntro}). This representation of $u^{\epsilon}(x)$ and $v^{\epsilon}(x)$ is unique up to multiplicative and additive constants. In fact one can multiply one of these functions by some
constant and divide the other function by the same constant or add a
constant to either of them.

Using the results in \cite{FW4} one can show (see Theorem 4.4 in \cite{Hyejin}) that
\begin{equation}
\hat{X}^{\epsilon}_{\cdot}\longrightarrow X_{\cdot}\textrm{ weakly in } \mathcal{C}_{0T}, \textrm{ for any } T<\infty, \textrm{ as }\epsilon\downarrow 0,\label{EscortProcessConvergence}
\end{equation}
where $X$ is the limiting process with operator defined by (\ref{LimitingProcess1}).
\vspace{0.5cm}

\begin{proof}[Proof of Lemma \ref{Lemma23}]
We prove the lemma just for $x\in[x_{1},x_{2}]$. Clearly, the proof for $x\in[-x_{2},-x_{1}]$ is the same.

We claim that it is sufficient to prove that
\begin{equation}
|\mathbb{E}^{\epsilon}_{x,y}[e^{-\lambda \tau^{\epsilon}(x_{1},x_{2})}f(X^{\epsilon}_{\tau^{\epsilon}(x_{1},x_{2})})+\int_{0}^{\tau^{\epsilon}(x_{1},x_{2})}e^{-\lambda t} [\lambda f(X^{\epsilon}_{t})-\hat{L}^{\epsilon}f(X^{\epsilon}_{t})]dt]- f(x)|\rightarrow 0,\label{Lemma23Alternative1}
\end{equation}
as $\epsilon \downarrow 0$, where $\hat{L}^{\epsilon}$ is defined in (\ref{OperatorOfHatProcess}). The left hand side of (\ref{Lemma23Alternative1}) is meaningful since $f$ is sufficiently smooth for $x\in[x_{1},x_{2}]$.

We observe that:
\begin{eqnarray}
& &|\mathbb{E}^{\epsilon}_{x,y}[e^{-\lambda \tau^{\epsilon}(x_{1},x_{2})}f(X^{\epsilon}_{\tau^{\epsilon}(x_{1},x_{2})})+\int_{0}^{\tau^{\epsilon}(x_{1},x_{2})}e^{-\lambda t} [\lambda f(X^{\epsilon}_{t})-Lf(X^{\epsilon}_{t})]dt]- f(x)|\nonumber\\
&\leq&|\mathbb{E}^{\epsilon}_{x,y}[e^{-\lambda \tau^{\epsilon}(x_{1},x_{2})}f(X^{\epsilon}_{\tau^{\epsilon}(x_{1},x_{2})})+\int_{0}^{\tau^{\epsilon}(x_{1},x_{2})}e^{-\lambda t} [\lambda f(X^{\epsilon}_{t})-\hat{L}^{\epsilon}f(X^{\epsilon}_{t})]dt]- f(x)|\nonumber\\
& & \hspace{0.2cm} + |\mathbb{E}^{\epsilon}_{x,y}\int_{0}^{\tau^{\epsilon}(x_{1},x_{2})}e^{-\lambda t}[Lf(X^{\epsilon}_{t})-\hat{L}^{\epsilon}f(X^{\epsilon}_{t})]dt|
\label{Lemma23Alternative2}
\end{eqnarray}
Now we claim that
\begin{equation}
\Vert \hat{L}^{\epsilon}f-Lf\Vert_{[x_{1},x_{2}]} \rightarrow 0, \textrm{ as } \epsilon\downarrow 0, \label{InequalityConvergenceOfGenerators}
\end{equation}
where for any function $g$ we define $\Vert g\Vert_{[x_{1},x_{2}]}=\sup_{x\in [x_{1},x_{2}]}|g(x)|$.
This follows directly by our assumptions on the function $V^{\epsilon}(x)$.
Therefore, it is indeed enough to prove (\ref{Lemma23Alternative1}).

By the It\^o formula applied to the function $e^{-\lambda t}f(x)$ we immediately get that (\ref{Lemma23Alternative1}) is equivalent to
\begin{equation}
|\mathbb{E}^{\epsilon}_{x,y}[\int_{0}^{\tau^{\epsilon}(x_{1},x_{2})}e^{-\lambda u}f_{x}(X^{\epsilon}_{u})\gamma_{1}^{\epsilon}(X^{\epsilon}_{u},Y^{\epsilon}_{u})dL^{\epsilon}_{u}-\int_{0}^{\tau^{\epsilon}(x_{1},x_{2})}\frac{1}{2}e^{-\lambda u} [f_{x}\frac{V^{\epsilon}_{x}}{V^{\epsilon}}](X_{u}^{\epsilon})du]|\rightarrow 0,\label{Lemma23Alternative3}
\end{equation}
as $\epsilon \downarrow 0$. We can estimate the left hand side of (\ref{Lemma23Alternative3}) as in Lemma 2.1 of \cite{FS1}:

Consider the auxiliary function
\begin{equation}
v^{\epsilon}(x,y)=\frac{1}{2}y^{2}f_{x}(x)\frac{V^{\epsilon}_{x}(x)}{V^{\epsilon}(x)}+yf_{x}(x)\frac{V^{u,\epsilon}_{x}(x)V^{l,\epsilon}(x)-V^{l,\epsilon}_{x}(x)V^{u,\epsilon}(x)}{V^{\epsilon}(x)}. \label{AuxillaryFunction2}
\end{equation}
It is easy to see that $v^{\epsilon}$ is a solution to the P.D.E.
\begin{eqnarray}
v^{\epsilon}_{yy}(x,y)&=&f_{x}(x)\frac{V^{\epsilon}_{x}(x)}{V^{\epsilon}(x)}, \hspace{0.2cm} y\in D^{\epsilon}_{x}\nonumber\\
\frac{\partial_{y} v^{\epsilon}(x,y)}{\partial n^{\epsilon}(x,y)}&=&-f_{x}(x)\frac{\gamma_{1}^{\epsilon}(x,y)}{|\gamma_{2}^{\epsilon}(x,y)|}, \hspace{0.2cm} y\in \partial D^{\epsilon}_{x}, \label{SpecificNeymmanProblem2}
\end{eqnarray}
where $n^{\epsilon}(x,y)=\frac{\gamma_{2}^{\epsilon}(x,y)}{|\gamma_{2}^{\epsilon}(x,y)|}$ and $ x \in \mathbb{R}$ is a parameter.

If we apply It\^{o} formula to the function $e^{-\lambda
t}v^{\epsilon}(x,y)$ we get that $e^{-\lambda t}
v^{\epsilon}(X_{t}^{\epsilon},Y_{t}^{\epsilon})$ satisfies with
probability one:
\begin{eqnarray}
e^{-\lambda t}v^{\epsilon}(X_{t}^{\epsilon},Y_{t}^{\epsilon})&=&v^{\epsilon}(x,y)+
\int_{0}^{t}e^{-\lambda s}[\frac{1}{2}\bigtriangleup v^{\epsilon}(X_{s}^{\epsilon},Y_{s}^{\epsilon})-\lambda v^{\epsilon}(X_{s}^{\epsilon},Y_{s}^{\epsilon})]ds\nonumber\\
&+&\int_{0}^{t}e^{-\lambda s}v^{\epsilon}_{x}(X_{s}^{\epsilon},Y_{s}^{\epsilon})dW_{s}^{1}
+\int_{0}^{t}e^{-\lambda s}v^{\epsilon}_{y}(X_{s}^{\epsilon},Y_{s}^{\epsilon})dW_{s}^{2}\nonumber\\
&+&\int_{0}^{t}e^{-\lambda s}v^{\epsilon}_{x}(X_{s}^{\epsilon},Y_{s}^{\epsilon})\gamma^{\epsilon}_{1}(X_{s}^{\epsilon},Y_{s}^{\epsilon})dL^{\epsilon}_{s}\nonumber\\
&+&\int_{0}^{t}e^{-\lambda s}v^{\epsilon}_{y}(X_{s}^{\epsilon},Y_{s}^{\epsilon})\gamma^{\epsilon}_{2}(X_{s}^{\epsilon},Y_{s}^{\epsilon})dL^{\epsilon}_{s}\label{ItoFormulaForLocalTime1}
\end{eqnarray}
Since $v^{\epsilon}(x,y)$ satisfies (\ref{SpecificNeymmanProblem2})
we have:
\begin{eqnarray}
& &
|\mathbb{E}^{\epsilon}_{x,y}[\int_{0}^{\tau^{\epsilon}(x_{1},x_{2})}e^{-\lambda
u}f_{x}(X^{\epsilon}_{u})\gamma_{1}^{\epsilon}(X^{\epsilon}_{u},Y^{\epsilon}_{u})dL^{\epsilon}_{u}-\int_{0}^{\tau^{\epsilon}(x_{1},x_{2})}\frac{1}{2}e^{-\lambda
u}
[f_{x}\frac{V^{\epsilon}_{x}}{V^{\epsilon}}](X_{u}^{\epsilon})du]|\leq \nonumber\\
&\leq& |\mathbb{E}^{\epsilon}_{x,y}e^{-\lambda
\tau^{\epsilon}(x_{1},x_{2})}v^{\epsilon}(X_{\tau^{\epsilon}(x_{1},x_{2})}^{\epsilon},Y_{\tau^{\epsilon}(x_{1},x_{2})}^{\epsilon})|+|v^{\epsilon}(x,y)|+\nonumber\\
&+&|\mathbb{E}^{\epsilon}_{x,y}[\int_{0}^{\tau^{\epsilon}(x_{1},x_{2})}e^{-\lambda s}[\frac{1}{2} v^{\epsilon}_{xx}(X_{s}^{\epsilon},Y_{s}^{\epsilon})-\lambda v^{\epsilon}(X_{s}^{\epsilon},Y_{s}^{\epsilon})]ds]|+\nonumber\\
&+&|\mathbb{E}^{\epsilon}_{x,y}[\int_{0}^{\tau^{\epsilon}(x_{1},x_{2})}e^{-\lambda
s}v^{\epsilon}_{x}(X_{s}^{\epsilon},Y_{s}^{\epsilon})\gamma^{\epsilon}_{1}(X_{s}^{\epsilon},Y_{s}^{\epsilon})dL^{\epsilon}_{s}]|\label{BoundByItoFormulaForLocalTime}
\end{eqnarray}
For any time $t\in[0,\tau^{\epsilon}(x_{1},x_{2})]$ the x-component
of the process $(X^{\epsilon}_{t},Y^{\epsilon}_{t})$ satisfies
$x_{1}\leq X^{\epsilon}_{t}\leq x_{2}$, i.e. it is far away from the
point of discontinuity. Taking into account the latter and the
definition of $v^{\epsilon}(x,y)$ by (\ref{AuxillaryFunction2}) we
get that the first three terms in the right hand side of
(\ref{BoundByItoFormulaForLocalTime}) are bounded by
$\epsilon^{2}C_{0}$ for $\epsilon$ small enough. So, it remains to
consider the last term, i.e. the integral in local time. First of
all it is easy to see that there exists a $C_{0}>0$ such that
\begin{equation}
|\mathbb{E}^{\epsilon}_{x,y}[\int_{0}^{\tau^{\epsilon}(x_{1},x_{2})}e^{-\lambda
s}v^{\epsilon}_{x}(X_{s}^{\epsilon},Y_{s}^{\epsilon})\gamma^{\epsilon}_{1}(X_{s}^{\epsilon},Y_{s}^{\epsilon})dL^{\epsilon}_{s}]|\leq
\epsilon^{2}C_{0}\mathbb{E}^{\epsilon}_{x,y}\int_{0}^{\tau^{\epsilon}(x_{1},x_{2})}
e^{-\lambda t} \epsilon
dL_{t}^{\epsilon}\label{BoundByItoFormulaForLocalTime1}
\end{equation}
As far as the integral in local time on the right hand side of
(\ref{BoundByItoFormulaForLocalTime1}) is concerned, we claim that
there exists an $\epsilon_{0}>0$ and a $C _{0}>0$ such that for all
$\epsilon<\epsilon_{0}$
\begin{equation}
\mathbb{E}^{\epsilon}_{x,y}\int_{0}^{\tau^{\epsilon}(x_{1},x_{2})}
e^{-\lambda t} \epsilon dL_{t}^{\epsilon}\leq
C_{0}.\label{BoundByItoFormulaForLocalTime2}
\end{equation}
Hence, taking into account (\ref{BoundByItoFormulaForLocalTime1})
and (\ref{BoundByItoFormulaForLocalTime2}) we get that the right
hand side of (\ref{BoundByItoFormulaForLocalTime}) converges to
zero. So the convergence (\ref{Lemma23Alternative3}) holds.

It remains to prove the claim
(\ref{BoundByItoFormulaForLocalTime2}). This can be done as in Lemma
2.2 of \cite{FS1}. In particular, one considers  the auxiliary
function
\begin{displaymath}
w^{\epsilon}(x,y)=y^{2}\frac{1}{V^{\epsilon}(x)}+y\frac{V^{l,\epsilon}(x)-V^{u,\epsilon}(x)}{V^{\epsilon}(x)}.
\end{displaymath}
It is easy to see that $w^{\epsilon}$ is a solution to the P.D.E.
\begin{eqnarray}
w^{\epsilon}_{yy}(x,y)&=&\frac{2}{V^{\epsilon}(x)}, \hspace{0.2cm} y\in D^{\epsilon}_{x}\nonumber\\
\frac{\partial_{y} w^{\epsilon}(x,y)}{\partial
n^{\epsilon}(x,y)}&=&-1, \hspace{0.2cm} y\in \partial
D^{\epsilon}_{x}, \label{SpecificNeymmanProblem2a}
\end{eqnarray}
where
$n^{\epsilon}(x,y)=\frac{\gamma_{2}^{\epsilon}(x,y)}{|\gamma_{2}^{\epsilon}(x,y)|}$
and $ x \in \mathbb{R}$ is a parameter.

The claim follows now directly by applying It\^{o} formula to the
function $e^{-\lambda t}w^{\epsilon}(x,y)$.
\end{proof}

\section{Proof of Lemma \ref{Lemma24}}

\begin{proof}[Proof of Lemma \ref{Lemma24}]

Let $\hat{\tau}^{\epsilon}=\hat{\tau}^{\epsilon}(\pm \kappa)$ be the
exit time of $\hat{X}^{\epsilon}_{t}$ (see (\ref{EscordProcess}))
from the interval $(-\kappa, \kappa)$. Denote also by
$\hat{\mathbb{E}}_{x}$ the mathematical expectation related to the
probability law induced by $\hat{X}^{\epsilon,x}_{t}$.

Let $\eta$ be a positive number, $\kappa>0$ be a small positive number and consider $x_{0}$ with $|x_{0}|< \kappa$ and $y_{0}\in D^{\epsilon}_{x_{0}}$. We want to estimate the difference
$$|\mathbb{E}^{\epsilon}_{x_{0},y_{0}}\tau^{\epsilon} (\pm \kappa)-\hat{\mathbb{E}}^{\epsilon}_{x_{0}}\hat{\tau}^{\epsilon}(\pm
\kappa)|.$$
As it can be derived by Theorem 2.5.1 of \cite{F1} the function
$\phi^{\epsilon}(x,y)=\mathbb{E}^{\epsilon}_{x,y}\tau^{\epsilon}(\pm
\kappa)$ is solution to
\begin{eqnarray}
\frac{1}{2}\bigtriangleup \phi^{\epsilon}(x,y)&=&-1 \hspace{0.5cm}\mathrm{ in }\hspace{0.2cm} (x,y)\in (-\kappa,\kappa)\times D^{\epsilon}_{x}\nonumber\\
\phi^{\epsilon}(\pm\kappa,y)&=& 0\nonumber\\
\frac{\partial \phi^{\epsilon}}{\partial \gamma^{\epsilon}}(x,y)&=& 0
\hspace{0.6cm}\mathrm{ on } \hspace{0.2cm}(x,y)\in (-\kappa,\kappa)\times \partial
D^{\epsilon}_{x}\label{PDEforLemma44_1}
\end{eqnarray}
Moreover,
$\hat{\phi}^{\epsilon}(x)=\hat{\mathbb{E}}^{\epsilon}_{x}\hat{\tau}^{\epsilon}(\pm
\kappa)$ is solution to
\begin{eqnarray}
\frac{1}{2} \hat{\phi}^{\epsilon}_{xx}(x)+\frac{1}{2}\frac{V_{x}^{\epsilon}(x)}{V^{\epsilon}(x)}\hat{\phi}^{\epsilon}_{x}(x)&=&-1 \hspace{0.5cm}\mathrm{ in }\hspace{0.2cm} x\in (-\kappa,\kappa)\nonumber\\
\hat{\phi}^{\epsilon}(\pm\kappa)&=& 0 \label{PDEforLemma44_2}
\end{eqnarray}
Let
$f^{\epsilon}(x,y)=\phi^{\epsilon}(x,y)-\hat{\phi}^{\epsilon}(x)$.
Then $f^{\epsilon}$ will satisfy
\begin{eqnarray}
\frac{1}{2}\bigtriangleup f^{\epsilon}(x,y)&=&\frac{1}{2}\frac{V_{x}^{\epsilon}(x)}{V^{\epsilon}(x)}\hat{\phi}^{\epsilon}_{x}(x) \hspace{0.5cm}\mathrm{ in }\hspace{0.2cm} (x,y)\in (-\kappa,\kappa)\times D^{\epsilon}_{x}\nonumber\\
f^{\epsilon}(\pm\kappa,y)&=&0\label{PDEforLemma44_3}\\
\frac{\partial f^{\epsilon}}{\partial \gamma^{\epsilon}}(x,y)&=& - \hat{\phi}^{\epsilon}_{x}(x)\gamma_{1}^{\epsilon}(x,y) \hspace{0.5cm}\mathrm{ on } \hspace{0.2cm}(x,y)\in (-\kappa,\kappa)\times \partial D^{\epsilon}_{x}\nonumber
\end{eqnarray}
By applying It\^o formula to the function $f^{\epsilon}$ and
recalling that $f^{\epsilon}$ satisfies (\ref{PDEforLemma44_3}) we
get that
\begin{equation}
|f^{\epsilon}(x,y)|\leq
|\mathbb{E}^{\epsilon}_{x,y}[\int_{0}^{\tau^{\epsilon}(\pm\kappa)}\hat{\phi}^{\epsilon}_{x}(X_{u}^{\epsilon})\gamma_{1}^{\epsilon}(X_{u}^{\epsilon},Y_{u}^{\epsilon})dL_{u}^{\epsilon}-\int_{0}^{\tau^{\epsilon}(\pm\kappa)}\frac{1}{2}[\hat{\phi}^{\epsilon}_{x}\frac{V^{\epsilon}_{x}}{V^{\epsilon}}](X_{u}^{\epsilon})du]|\label{PDEforLemma44_4}
\end{equation}
We can estimate the right hand side of (\ref{PDEforLemma44_4})
similarly to the left hand side of (\ref{Lemma23Alternative3}) of
Lemma \ref{Lemma23} (see also Lemma 2.1 in \cite{FS1}). Consider the
auxiliary function
\begin{equation}
w^{\epsilon}(x,y)=\frac{1}{2}y^{2}\hat{\phi}^{\epsilon}_{x}(x)\frac{V^{\epsilon}_{x}(x)}{V^{\epsilon}(x)}+y\hat{\phi}^{\epsilon}_{x}(x)\frac{V^{u,\epsilon}_{x}(x)V^{l,\epsilon}(x)-V^{l,\epsilon}_{x}(x)V^{u,\epsilon}(x)}{V^{\epsilon}(x)} \label{AuxillaryFunction1}
\end{equation}
It is easy to see that $w$ is a solution to the P.D.E.
\begin{eqnarray}
w^{\epsilon}_{yy}(x,y)&=&\hat{\phi}^{\epsilon}_{x}(x)\frac{V^{\epsilon}_{x}(x)}{V^{\epsilon}(x)}, \hspace{0.2cm} y\in D^{\epsilon}_{x}\nonumber\\
\frac{\partial_{y} w^{\epsilon}(x,y)}{\partial n^{\epsilon}(x,y)}&=&-\hat{\phi}^{\epsilon}_{x}(x)\frac{\gamma_{1}^{\epsilon}(x,y)}{|\gamma_{2}^{\epsilon}(x,y)|}, \hspace{0.2cm} y\in \partial D^{\epsilon}_{x}, \label{SpecificNeymmanProblem3}
\end{eqnarray}
where $n^{\epsilon}(x,y)=\frac{\gamma_{2}^{\epsilon}(x,y)}{|\gamma_{2}^{\epsilon}(x,y)|}$ and $ x \in \mathbb{R}$ is a parameter.

Then if we apply It\^{o} formula to the function $w^{\epsilon}(x,y)$
and recall that $w^{\epsilon}$ satisfies
(\ref{SpecificNeymmanProblem3}), we get an upper bound for the right
hand side of (\ref{PDEforLemma44_4}) that is the same to the right
hand side of (\ref{BoundByItoFormulaForLocalTime}) with $\lambda=0$,
$v^{\epsilon}$ replaced by $w^{\epsilon}$ and
$\tau^{\epsilon}(x_{1},x_{2})$ replaced by $\tau^{\epsilon}(\pm
\kappa)$. Namely, for $(x,y)=(x_{0},y_{0})$, we have the following
\begin{eqnarray}
& &
|\mathbb{E}^{\epsilon}_{x_{0},y_{0}}[\int_{0}^{\tau^{\epsilon}(\pm
\kappa)}\hat{\phi}^{\epsilon}_{x}(X_{u}^{\epsilon})\gamma_{1}^{\epsilon}(X_{u}^{\epsilon},Y_{u}^{\epsilon})dL_{u}^{\epsilon}-\int_{0}^{\tau^{\epsilon}(\pm
\kappa)}\frac{1}{2}[\hat{\phi}^{\epsilon}_{x}\frac{V^{\epsilon}_{x}}{V^{\epsilon}}](X_{u}^{\epsilon})du]|\leq \nonumber\\
&\leq& \sup_{(x,y)\in\lbrace\pm\kappa\rbrace\times D^{\epsilon}_{\pm \kappa}}|w^{\epsilon}(x,y)|+\sup_{(x,y)\in\lbrace x_{0}\rbrace\times D^{\epsilon}_{ x_{0}}}|w^{\epsilon}(x,y)|+\label{BoundByItoFormulaForLocalTimeLemma24}\\
&+&|\mathbb{E}^{\epsilon}_{x_{0},y_{0}}\int_{0}^{\tau^{\epsilon}(\pm
\kappa)}\frac{1}{2} w^{\epsilon}_{xx}(X_{s}^{\epsilon},Y_{s}^{\epsilon})ds|+|\mathbb{E}^{\epsilon}_{x_{0},y_{0}}\int_{0}^{\tau^{\epsilon}(\pm
\kappa)} w^{\epsilon}_{x}(X_{s}^{\epsilon},Y_{s}^{\epsilon})\gamma^{\epsilon}_{1}(X_{s}^{\epsilon},Y_{s}^{\epsilon})dL^{\epsilon}_{s}|\nonumber
\end{eqnarray}
Now one can solve (\ref{PDEforLemma44_2}) explicitly and get that for $x\in[-\kappa,\kappa]$
\begin{equation}
\hat{\phi}^{\epsilon}(x)=\int_{-\kappa}^{x}\frac{-2}{V^{\epsilon}(y)}\int_{-\kappa}^{y}V^{\epsilon}(z)dzdy+[\frac{\int_{-\kappa}^{\kappa}\frac{2}{V^{\epsilon}(y)}
\int_{-\kappa}^{y}V^{\epsilon}(z)dzdy}{\int_{-\kappa}^{\kappa}\frac{1}{V^{\epsilon}(y)}dy}]\int_{-\kappa}^{x}\frac{1}{V^{\epsilon}(y)}dy.\label{TheHatExitTime}
\end{equation}
Using (\ref{TheHatExitTime}) and the form of $V^{\epsilon}(x)$ as
described by
(\ref{DefinitionOfCrossSections1})-(\ref{DefinitionOfCrossSections5})
we get that the first two terms of the right hand side of
(\ref{BoundByItoFormulaForLocalTimeLemma24}) can be made arbitrarily small  for $\epsilon$ sufficiently small. For $\epsilon$ small enough the two integral terms of
the right hand side of (\ref{BoundByItoFormulaForLocalTimeLemma24})
are bounded by $C_{0}\xi^{\epsilon}
\mathbb{E}^{\epsilon}_{x_{0},y_{0}}\tau^{\epsilon}(\pm \kappa)$, where $\xi^{\epsilon}$ is defined in (\ref{DefinitionOfCrossSections5}). The
local time integral can be treated as in Lemma 2.2 of \cite{FS1}, so
it will not be repeated here (see also the end of the proof of Lemma
\ref{Lemma23}). In reference to the latter, we mention that the singularity at the point $x=0$ complicates a bit the situation. However, assumption (\ref{DefinitionOfCrossSections5}) allows one to follow the procedure mentioned and derive the aforementioned estimate for the local time integral.  Hence, we have the following upper bound for
$f^{\epsilon}(x_{0},y_{0})$
\begin{equation}
|f^{\epsilon}(x_{0},y_{0})|\leq C_{0}[\kappa \eta
+\xi^{\epsilon}
\mathbb{E}^{\epsilon}_{x_{0},y_{0}}\tau^{\epsilon}(\pm
\kappa)].\label{DifferenceForHatAndXprocess}
\end{equation}
Moreover, it follows from (\ref{TheHatExitTime}) (see also \cite{Hyejin}) that for $\eta>0$ there exists a $\kappa_{\eta}>0$ such that for every $0<\kappa<\kappa_{\eta}$, for sufficiently small $\epsilon$ and for all $x$ with $|x|\leq \kappa$
\begin{equation}
\frac{1}{\kappa}|\hat{\mathbb{E}}^{\epsilon}_{
x}\hat{\tau}^{\epsilon}(\pm \kappa)-\kappa \theta|\leq \eta
\label{ExitTimeForEscordProcessWithDelay}
\end{equation}
,where $\theta=\frac{v(0+)-v(0-)}{[u'(0+)]^{-1}+[u'(0-)]^{-1}}$.

Therefore, since $\xi^{\epsilon}\downarrow 0$ (by
assumption (\ref{DefinitionOfCrossSections5})),
(\ref{DifferenceForHatAndXprocess}) and
(\ref{ExitTimeForEscordProcessWithDelay}) give us for sufficiently small $\epsilon$ that
\begin{displaymath}
\frac{1}{\kappa}|\mathbb{E}^{\epsilon}_{
x_{0},y_{0}}\tau^{\epsilon}(\pm
\kappa)-\hat{\mathbb{E}}^{\epsilon}_{
x_{0}}\hat{\tau}^{\epsilon}(\pm \kappa)|\leq C_{0}\eta.
\end{displaymath}
The latter and (\ref{ExitTimeForEscordProcessWithDelay}) conclude the proof of the lemma.
\end{proof}

\section{Proof of Lemma \ref{Lemma25}}

In order to prove Lemma \ref{Lemma25} we will make use of a result regarding the
invariant measures of the associated Markov chains (see Lemma
\ref{Lemma52}) and a result regarding the strong Markov character of
the limiting process (see Lemma \ref{Lemma54} and the beginning of
the proof of Lemma \ref{Lemma25}).

Of course, since the gluing conditions at $0$ are of local
character, it is sufficient to consider not the whole domain
$D^{\epsilon}$, but just the part of $D^{\epsilon}$ that is in the
neighborhood of $x=0$. Thus, we consider the process
$(X^{\epsilon}_{t},Y^{\epsilon}_{t})$ in
\begin{displaymath}
\Xi^{\epsilon}=\lbrace (x,y):|x|\leq 1, y \in D^{\epsilon}_{x}\rbrace
\end{displaymath}
that reflects normally on $\partial \Xi^{\epsilon}$.

Recall that  $0<\kappa_{0}<\kappa$. Define the set
\begin{equation}
\Gamma_{\kappa}=\lbrace (x,y):x=\kappa, y \in
D^{\epsilon}_{\kappa}\rbrace \label{DefinitionOfGammaSet}
\end{equation}
For notational convenience we will write $\Gamma_{\pm \kappa}=\Gamma_{\kappa}\cup \Gamma_{-\kappa}$.

Define $\Delta_{\kappa}=\Gamma_{\pm\kappa}\cup \Gamma_{\pm(1-\kappa)}$ and $\Delta_{\kappa_{0}}=\Gamma_{\pm\kappa_{0}}\cup \Gamma_{\pm(1-\kappa_{0})}$. We consider two cycles of Markov times $\lbrace \tau_{n}\rbrace$ and $\lbrace
\sigma_{n}\rbrace$ such that:
\begin{displaymath}
0=\sigma_{0}\leq\tau_{0}<\sigma_{1}<\tau_{1}<\sigma_{2}<\tau_{2}<\dots
\end{displaymath}
where:
\begin{eqnarray}
\tau_{n}&=&\inf \lbrace t\geq\sigma_{n}: (X^{\epsilon}_{t}, Y_{t}^{\epsilon}) \in \Delta_{\kappa} \rbrace \nonumber\\
\sigma_{n}&=&\inf \lbrace t\geq\tau_{n-1}: (X^{\epsilon}_{t},
Y_{t}^{\epsilon}) \in \Delta_{\kappa_{0}} \rbrace
\label{MarkovTimes3}
\end{eqnarray}
We will use the relations
\begin{eqnarray}
\mu_{\epsilon}(A)&=&\int_{\Delta_{\kappa}}\nu_{\epsilon}^{\kappa}(dx,dy)\mathbb{E}^{\epsilon}_{x,y}\int_{0}^{\tau_{1}}
\chi_{[(X^{\epsilon}_{t},Y^{\epsilon}_{t})\in
A]}dt\nonumber\\
&=&\int_{\Delta_{\kappa_{0}}}\nu_{\epsilon}^{\kappa_{0}}(dx,dy)\mathbb{E}^{\epsilon}_{x,y}\int_{0}^{\sigma_{1}}
\chi_{[(X^{\epsilon}_{t},Y^{\epsilon}_{t})\in
A]}dt\label{InvariantMeasures}
\end{eqnarray}
and
\begin{eqnarray}
\nu_{\epsilon}^{\kappa}(B)&=&\int_{\Delta_{\kappa_{0}}}\nu_{\epsilon}^{\kappa_{0}}(dx,dy)\mathbb{P}^{\epsilon}_{x,y}[(X^{\epsilon}_{\tau_{0}},Y^{\epsilon}_{\tau_{0}})\in
B]\nonumber\\
\nu_{\epsilon}^{\kappa_{0}}(C)&=&\int_{\Delta_{\kappa}}\nu_{\epsilon}^{\kappa}(dx,dy)\mathbb{P}^{\epsilon}_{x,y}[(X^{\epsilon}_{\sigma_{1}},Y^{\epsilon}_{\sigma_{1}})\in
C]\label{InvariantMeasuresForMarkovChains}
\end{eqnarray}
between the invariant measures $\mu_{\epsilon}$ of the process
$(X^{\epsilon}_{t},Y^{\epsilon}_{t})$ in $D^{\epsilon}$,
$\nu_{\epsilon}^{\kappa}$ and $\nu_{\epsilon}^{\kappa_{0}}$ of the
Markov chains $(X^{\epsilon}_{\tau_{n}},Y^{\epsilon}_{\tau_{n}})$
and $(X^{\epsilon}_{\sigma_{n}},Y^{\epsilon}_{\sigma_{n}})$ on
$\Delta_{\kappa}$ and $\Delta_{\kappa_{0}}$ respectively (see
\cite{K1}; also \cite{FW1}). It is clear that the first invariant
measure is, up to a constant, the Lebesgue measure on
$D^{\epsilon}$.

Formula (\ref{InvariantMeasures}) implies the corresponding equality
for the integrals with respect to the invariant measure
$\mu_{\epsilon}$. This means that if $\Psi(x,y)$ is an integrable
function, then
\begin{eqnarray}
\int\int_{\mathbb{R}^{2}} \Psi(x,y)\mu_{\epsilon}(dx,dy)&=&\int\int_{\Xi^{\epsilon}} \Psi(x,y)dxdy\nonumber\\
&=&\int_{\Delta_{\kappa}}\nu_{\epsilon}^{\kappa}(dx,dy)\mathbb{E}^{\epsilon}_{x,y}\int_{0}^{\tau_{1}}
\Psi(X^{\epsilon}_{t},Y^{\epsilon}_{t})dt\nonumber\\
&=&\int_{\Delta_{\kappa_{0}}}\nu_{\epsilon}^{\kappa_{0}}(dx,dy)\mathbb{E}^{\epsilon}_{x,y}\int_{0}^{\sigma_{1}}
\Psi(X^{\epsilon}_{t},Y^{\epsilon}_{t})dt\label{InvariantMeasuresWithPsi}
\end{eqnarray}
The following simple lemma will be used
in the proof of Lemmata \ref{Lemma52} and \ref{Lemma54}.

\begin{lem}
Let $0<x_{1}<x_{2}$, $\psi$ be a function defined in $[x_{1},x_{2}]$ and $\phi$ be a function defined on $x_{1}$ and $x_{2}$. Then
\begin{displaymath}
\lim_{\epsilon\downarrow 0} \mathbb{E}^{\epsilon}_{x,y}[\phi(X^{\epsilon}_{\tau(x_{1},x_{2})})+\int_{0}^{\tau(x_{1},x_{2})}\psi(X^{\epsilon}_{t})dt]=g(x),
\end{displaymath}
uniformly in $(x,y)\in [x_{1},x_{2}]\times D^{\epsilon}_{x}$  and
\begin{eqnarray}
g(x)&=&\frac{u(x_{2})-u(x)}{u(x_{2})-u(x_{1})}[\phi(x_{1})+\int_{x_{1}}^{x}(u(y)-u(x_{1}))\psi(y)dv(y)]\nonumber\\
&+&\frac{u(x)-u(x_{1})}{u(x_{2})-u(x_{1})}[\phi(x_{2})+\int_{x}^{x_{2}}(u(x_{2})-u(y))\psi(y)dv(y)]\nonumber
\end{eqnarray}
A similar result holds for $(x,y)\in [-x_{2},-x_{1}]\times D^{\epsilon}_{x}$. \label{Lemma51}
\end{lem}
\begin{proof}
This lemma is similar to Lemma 8.4.6 of \cite{FW1}, so we briefly
outline its proof. First one proves that
$\mathbb{E}^{\epsilon}_{x,y}\tau(x_{1},x_{2})$ is bounded in
$\epsilon$ for $\epsilon$ small enough and for all $(x,y)\in[x_{1},x_{2}]\times D^{\epsilon}_{x}$. The
latter and the proof of Lemma \ref{Lemma23} show that in this case we can take $\lambda=0$ in Lemma (\ref{Lemma23}) and apply it to the function $g$ that
is the solution to
\begin{eqnarray}
D_{v}D_{u}g(x)&=&-\psi(x), \hspace{0.2cm} x_{1}<x<x_{2}\nonumber\\
g(x_{i})&=&\phi(x_{i}), \hspace{0.2cm} i=1,2.\nonumber
\end{eqnarray}
This gives the desired result.
\end{proof}
Lemma \ref{Lemma52} below characterizes the asymptotics of the invariant measures
 $\nu_{\epsilon}^{\kappa}$ and $\nu_{\epsilon}^{\kappa_{0}}$.

\begin{lem}
Let $v$ be the function defined by (\ref{uANDvForLimitingProcessIntro}).  The following statements hold
\begin{eqnarray}
\lim_{\epsilon \downarrow 0} \frac{1}{\epsilon}\nu_{\epsilon}^{\kappa}( \Gamma_{\kappa})[u(\kappa)-u(\kappa_{0})]&=&1\nonumber\\
\lim_{\epsilon \downarrow 0}
\frac{1}{\epsilon}\nu_{\epsilon}^{\kappa}(
\Gamma_{1-\kappa})[u(1-\kappa_{0})-u(1-\kappa)]&=&1\label{AsymptoticsOfInvariantMeasures1}
\end{eqnarray}
Similar statements are true for  $\nu_{\epsilon}^{\kappa}(
\Gamma_{-\kappa}), \nu_{\epsilon}^{\kappa}(\Gamma_{-(1-\kappa)})$
and $\nu_{\epsilon}^{\kappa_{0}}$. \label{Lemma52}
\end{lem}

\begin{proof}
We calculate the asymptotics of the invariant measures
$\nu_{\epsilon}^{\kappa}$ and $\nu_{\epsilon}^{\kappa_{0}}$ by
selecting $\Psi$ properly. For this purpose let us choose
$\Psi(x,y)=\Psi(x)$ (i.e. it is a function only of $x$) that is
bounded, continuous and is $0$ outside $\kappa<x<1-\kappa$.

With this choice for $\Psi$ the left hand side of (\ref{InvariantMeasuresWithPsi}) becomes:
\begin{equation}
\int\int_{\Xi^{\epsilon}} \Psi(x,y)dxdy=\int_{\kappa}^{1-\kappa}\Psi(x)V^{\epsilon}(x)dx=
\epsilon\int_{\kappa}^{1-\kappa}\Psi(x)dv^{\epsilon}(x), \label{LHSWithPsi1}
\end{equation}
where we recall that $v^{\epsilon}(x)=\int_{0}^{x}\frac{V^{\epsilon}(y)}{\epsilon}dy$ (see (\ref{uANDvForEscortProcessIntro})).

Moreover, the particular choice if $\Psi$, also implies that
$\int_{\sigma_{1}}^{\tau_{1}} \Psi(X^{\epsilon}_{t})dt=0$. Then, the right hand side of (\ref{InvariantMeasuresWithPsi}) becomes:
\begin{eqnarray}
&
&\int_{\Delta_{\kappa}}\nu_{\epsilon}^{\kappa}(dx,dy)\mathbb{E}^{\epsilon}_{x,y}\int_{0}^{\tau_{1}}
\Psi(X^{\epsilon}_{t})dt=\int_{\Delta_{\kappa}}\nu_{\epsilon}^{\kappa}(dx,dy)\mathbb{E}^{\epsilon}_{x,y}\int_{0}^{\sigma_{1}}
\Psi(X^{\epsilon}_{t})dt\label{RHSWithPsi1}\\
&=&\int_{D^{\epsilon}_{\kappa}}\nu_{\epsilon}^{\kappa}(\kappa,dy)\mathbb{E}^{\epsilon}_{
\kappa,y}\int_{0}^{\sigma_{1}}
\Psi(X^{\epsilon}_{t})dt+\int_{D^{\epsilon}_{1-\kappa}}\nu_{\epsilon}^{\kappa}(1-
\kappa,dy)\mathbb{E}^{\epsilon}_{ 1-\kappa,y}\int_{0}^{\sigma_{1}}
\Psi(X^{\epsilon}_{t})dt\nonumber
\end{eqnarray}
Next, we express the right hand side of (\ref{RHSWithPsi1}) through
the $v$ and $u$ functions (defined by
(\ref{uANDvForLimitingProcessIntro})) using Lemma \ref{Lemma51}. We start
with the term $\mathbb{E}^{\epsilon}_{
\kappa,y}\int_{0}^{\sigma_{1}} \Psi(X^{\epsilon}_{t})dt$.  Use Lemma
\ref{Lemma51} with $\phi=0$ and $\psi(x)=\Psi(x)$. For sufficiently
small $\epsilon$ we have
\begin{eqnarray}
& &\mathbb{E}^{\epsilon}_{ \kappa,y}\int_{0}^{\sigma_{1}}
\Psi(X^{\epsilon}_{t})dt=
\frac{u(1-\kappa_{0})-u(\kappa)}{u(1-\kappa_{0})-u(\kappa_{0})}\int_{\kappa_{0}}^{\kappa}(u(y)-u(
\kappa_{0}))\Psi(y)dv(y)\nonumber\\
& &+ \frac{u(\kappa)-u(\kappa_{0})}{u(1-\kappa_{0})-u(\kappa_{0})}\int_{\kappa}^{1-\kappa_{0}}(u(1-\kappa_{0})-u(y))\Psi(y)dv(y)+o(1)\nonumber\\
& &=\frac{u(\kappa)-u(\kappa_{0})}{u(1-\kappa_{0})-u(\kappa_{0})}\int_{\kappa}^{1-\kappa}(u(1-\kappa_{0})-u(y))\Psi(y)dv(y)+o(1)\label{RHSWithPsi2}
\end{eqnarray}
where the term $o(1)\downarrow 0$ as $\epsilon \downarrow 0$.
Similarly for the term
$\mathbb{E}^{\epsilon}_{1-\kappa,y}\int_{0}^{\sigma_{1}}
\Psi(X^{\epsilon}_{t})dt$ we have
\begin{eqnarray}
& &\mathbb{E}^{\epsilon}_{ 1-\kappa,y}\int_{0}^{\sigma_{1}} \Psi(X^{\epsilon}_{t})dt=\frac{u(1-\kappa_{0})-u(1-\kappa)}{u(1-\kappa_{0})-u(\kappa_{0})}\int_{\kappa_{0}}^{1-\kappa}(u(y)-u(\kappa_{0}))\Psi(y)dv(y)\nonumber\\
& &+ \frac{u(1-\kappa)-u(\kappa_{0})}{u(1-\kappa_{0})-u(\kappa_{0})}\int_{1-\kappa}^{1-\kappa_{0}}(u(1-\kappa_{0})-u(y))\Psi(y)dv(y)+o(1)\nonumber\\
& &=\frac{u(1-\kappa_{0})-u(1-\kappa)}{u(1-\kappa_{0})-u(\kappa_{0})}\int_{\kappa}^{1-\kappa}(u(y)-u(\kappa_{0}))\Psi(y)dv(y)+o(1)\label{RHSWithPsi3}
\end{eqnarray}
Taking into account relations (\ref{InvariantMeasuresWithPsi}) and (\ref{LHSWithPsi1})-(\ref{RHSWithPsi3}) and the particular choice of the function $\Psi$ we get the following relation for sufficiently small $\epsilon$
\begin{eqnarray}
& & \epsilon\int_{\kappa}^{1-\kappa}\Psi(x)dv^{\epsilon}(x)=\label{InvariantMeasuresWithPsiSpecific1}\\
&=&\nu_{\epsilon}^{\kappa}( \Gamma_{\kappa})
\left[\frac{u(\kappa)-u(\kappa_{0})}{u(1-\kappa_{0})-u(\kappa_{0})}\int_{\kappa}^{1-\kappa}(u(1-\kappa_{0})-u(y))\Psi(y)dv(y)+o(1)\right]\nonumber\\
&+&\nu_{\epsilon}^{\kappa}(\Gamma_{1-\kappa})
\left[\frac{u(1-\kappa_{0})-u(1-\kappa)}{u(1-\kappa_{0})-u(\kappa_{0})}\int_{\kappa}^{1-\kappa}(u(y)-u(\kappa_{0}))\Psi(y)dv(y)+o(1)\right]\nonumber
\end{eqnarray}
At each continuity point of  $v(x)$ we have  $v(x)=\lim_{\epsilon \downarrow 0 }v^{\epsilon}(x)$, so the equality above is true for an arbitrary continuous function $\Psi$ if the following hold:
\begin{eqnarray}
\lim_{\epsilon \downarrow 0} \frac{1}{\epsilon}\nu_{\epsilon}^{\kappa}(\Gamma_{ \kappa})[u(\kappa)-u(\kappa_{0})]&=&1\nonumber\\
\lim_{\epsilon \downarrow 0}
\frac{1}{\epsilon}\nu_{\epsilon}^{\kappa}(\Gamma_{
1-\kappa})[u(1-\kappa_{0})-u(1-\kappa)]&=&1\nonumber
\end{eqnarray}
Thus, the proof of Lemma \ref{Lemma52} is complete.
\end{proof}
\begin{rem}
Equalities (\ref{InvariantMeasuresForMarkovChains}) immediately give
us that
$\nu_{\epsilon}^{\kappa}(\Gamma_{\pm\kappa})=\nu_{\epsilon}^{\kappa_{0}}(\Gamma_{\pm\kappa_{0}})$
and
$\nu_{\epsilon}^{\kappa}(\Gamma_{\pm(1-\kappa)})=\nu_{\epsilon}^{\kappa_{0}}(\Gamma_{\pm(1-\kappa_{0})})$
(see also the related discussion in chapter $8$ of \cite{FW1}).
Finally, it is easy to see that
$\nu_{\epsilon}^{\kappa}(\Gamma_{\pm\kappa})=\nu_{\epsilon}^{\kappa}(\Gamma_{\kappa})+\nu_{\epsilon}^{\kappa}(\Gamma_{-\kappa})$.\label{Remark43}
\end{rem}
\begin{flushright}
$\square$
\end{flushright}

\begin{lem} Let us consider fixed numbers
$0<x_{1}<x_{2}$ and let $x_{0}\in (x_{1},x_{2})$. For every $x_{0}$ we have
\begin{displaymath}
\lim_{\epsilon\downarrow 0}\max_{f:\parallel f\parallel\leq 1}
|\mathbb{E}^{\epsilon}_{x_{0},y_{01}}f(X^{\epsilon}_{\tau(x_{1},x_{2})},Y^{\epsilon}_{\tau(x_{1},x_{2})})-\mathbb{E}^{\epsilon}_{x_{0},y_{02}}f(X^{\epsilon}_{\tau(x_{1},x_{2})},Y^{\epsilon}_{\tau(x_{1},x_{2})})|=0,
\end{displaymath}
uniformly in $y_{01},y_{02}\in D^{\epsilon}_{x_{0}}$ and for functions $f(x,y)$ that are well defined on $(x,y)\in\{x_{1},x_{2}\}\times D^{\epsilon}_{x}$.
\label{Lemma54}
\end{lem}
\begin{proof}
We only need to observe that $(a)$:  Lemma \ref{Lemma51} applied to $\phi(x)=f(x,0)$ and $\psi(x)=0$ immediately gives us that
\begin{displaymath}
\lim_{\epsilon\downarrow 0}\max_{f:\parallel f\parallel\leq 1}
|\mathbb{E}^{\epsilon}_{x_{0},y_{01}}f(X^{\epsilon}_{\tau(x_{1},x_{2})},0)-\mathbb{E}^{\epsilon}_{x_{0},y_{02}}f(X^{\epsilon}_{\tau(x_{1},x_{2})},0)|=0,
\end{displaymath}
and $(b)$: the $y-$component of the process converges to zero.
\end{proof}

\begin{proof}[Proof of Lemma \ref{Lemma25}.] Firstly, we prove (using Lemma \ref{Lemma54}) that
$\mathbb{P}^{\epsilon}_{x,y}(X^{\epsilon}_{\tau(\pm \kappa)}
=\kappa)$ has approximately the same value for all $x$ that belong
to a small neighborhood of zero for $\epsilon$ small enough.
Secondly, we identify (using (\ref{InvariantMeasures}) and Lemma
\ref{Lemma52}) that this value is
$\frac{[u'(0+)]^{-1}}{[u'(0+)]^{-1}+[u'(0-)]^{-1}}$.

\vspace{0.3cm}

We begin by showing that for any $(x_{i},y_{i})\in \Gamma_{\pm \kappa_{0}}$ with $i=1,2$
\begin{equation}
|\mathbb{E}^{\epsilon}_{x_{1},y_{1}}\chi_{(X^{\epsilon}_{\tau(\pm
\kappa)}=\kappa)}-\mathbb{E}^{\epsilon}_{x_{2},y_{2}}\chi_{(X^{\epsilon}_{\tau(\pm
\kappa)}=\kappa)}|\rightarrow 0 \textrm{ as } \epsilon\downarrow 0,
\label{MaxMinProbAlternativeStatement}
\end{equation}
uniformly in $y_{1},y_{2}$.

Let us define for notational convenience
\begin{displaymath}
F^{\epsilon}(x,y)=\mathbb{E}^{\epsilon}_{x,y}[\chi_{(X^{\epsilon}_{\tau(\pm
\kappa)}=\kappa)}]
\end{displaymath}
Firstly, we prove that (\ref{MaxMinProbAlternativeStatement}) holds for
$(\kappa_{0},y_{1}),(\kappa_{0},y_{2})\in \Gamma_{\kappa_{0}}$. Let
$\kappa_{1},\kappa_{2}$ be such that
$0<\kappa_{1}<\kappa_{0}<\kappa_{2}<\kappa$.  Using strong
Markov property with respect to $\tau( \kappa_{1}, \kappa_{2})<
\tau(\pm \kappa)$ we get
\begin{equation}
\mathbb{E}^{\epsilon}_{\kappa_{0},y}\chi_{(X^{\epsilon}_{\tau(\pm
\kappa)}=\kappa)}=
\mathbb{E}^{\epsilon}_{\kappa_{0},y}\mathbb{E}^{\epsilon}_{X_{\tau(\kappa_{1},\kappa_{2})},Y_{\tau(\kappa_{1},\kappa_{2})}}[\chi_{(X^{\epsilon}_{\tau(\pm
\kappa)}=\kappa)}]\label{StrongMarkovPropertyForMaxMinProb}
\end{equation}
Equation (\ref{StrongMarkovPropertyForMaxMinProb}) and
Lemma \ref{Lemma54} imply that
(\ref{MaxMinProbAlternativeStatement}) holds for
$(\kappa_{0},y_{1}),(\kappa_{0},y_{2})\in \Gamma_{\kappa_{0}}$ uniformly in $y_{1},y_{2}\in D_{\kappa_{0}}^{\epsilon}$.

Hence, we have
\begin{equation}
|F^{\epsilon}(\kappa_{0},y_{1})-F^{\epsilon}(\kappa_{0},y_{2})|\rightarrow 0 \textrm{ as } \epsilon\downarrow 0.
\label{MaxMinProbAlternativeStatement1a}
\end{equation}
Similarly, it can be shown that (\ref{MaxMinProbAlternativeStatement1a}) holds for $(-\kappa_{0},y_{1}),(-\kappa_{0},y_{2})\in
\Gamma_{-\kappa_{0}}$ as well.

Secondly, we observe that if (\ref{MaxMinProbAlternativeStatement})
holds for $(x_{1},y_{1})=(\kappa_{0},0)$ and
$(x_{2},y_{2})=(-\kappa_{0},0)$ then it will hold for any
$(x_{i},y_{i})\in \Gamma_{\pm \kappa_{0}}$ with $i=1,2$. Indeed, we
have
\begin{eqnarray}
|F^{\epsilon}(\kappa_{0},y_{1})-F^{\epsilon}(-\kappa_{0},y_{2})|&\leq&|F^{\epsilon}(\kappa_{0},0)-F^{\epsilon}(-\kappa_{0},0)|+|F^{\epsilon}(\kappa_{0},y_{1})-F^{\epsilon}(\kappa_{0},0)|\nonumber\\
& &+|F^{\epsilon}(-\kappa_{0},y_{2})-F^{\epsilon}(-\kappa_{0},0)|
\nonumber
\end{eqnarray}
The last two terms in the right hand side of the inequality above converge to zero as $\epsilon\downarrow 0$ by the discussion above. Hence, it remains to prove that
\begin{equation}
|F^{\epsilon}(\kappa_{0},0)-F^{\epsilon}(-\kappa_{0},0)|\rightarrow
0 \textrm{ as } \epsilon\downarrow 0.
\label{MaxMinProbAlternativeStatement1}
\end{equation}
Let us choose some $\kappa^{'}_{0}$ such that
$0<\kappa_{0}<\kappa^{'}_{0}<\kappa$.
Obviously, if the process starts from some point on
$\Gamma_{x}$ with $x\in[-\kappa, \kappa_{0}]$ then
$\tau(-\kappa,\kappa^{'}_{0})\leq\tau(\pm\kappa)$. So, by applying
the strong Markov property with respect to
$\tau(-\kappa,\kappa^{'}_{0})\leq\tau(\pm\kappa)$ we have
\begin{eqnarray}
\inf_{y\in D^{\epsilon}_{\kappa^{'}_{0}}} F^{\epsilon}(\kappa^{'}_{0},y) \mathbb{P}^{\epsilon}_{x,0}[(X^{\epsilon}_{\tau(-\kappa,\kappa^{'}_{0})},Y^{\epsilon}_{\tau(-\kappa,\kappa^{'}_{0})}) \in \Gamma_{\kappa^{'}_{0}}]&\leq& F^{\epsilon}(x,0)\leq\nonumber\\
\leq\sup_{y\in D^{\epsilon}_{\kappa^{'}_{0}}} F^{\epsilon}(\kappa^{'}_{0},y)  \mathbb{P}^{\epsilon}_{x,0}[(X^{\epsilon}_{\tau(-\kappa,\kappa^{'}_{0})},Y^{\epsilon}_{\tau(-\kappa,\kappa^{'}_{0})}) \in \Gamma_{\kappa^{'}_{0}}]& &\label{MaxMinProbAlternativeStatement2b}
\end{eqnarray}
Using the latter, we have
\begin{eqnarray}
& &|F^{\epsilon}(\kappa_{0},0)-F^{\epsilon}(-\kappa_{0},0)|\leq \nonumber\\
&\leq&|\sup_{y\in D^{\epsilon}_{\kappa^{'}_{0}}} F^{\epsilon}(\kappa^{'}_{0},y)-\inf_{y\in D^{\epsilon}_{\kappa^{'}_{0}}} F^{\epsilon}(\kappa^{'}_{0},y)|\max_{x=\kappa_{0},-\kappa_{0}}\{\mathbb{P}^{\epsilon}_{x,0}[(X^{\epsilon}_{\tau(-\kappa,\kappa^{'}_{0})},Y^{\epsilon}_{\tau(-\kappa,\kappa^{'}_{0})}) \in \Gamma_{\kappa^{'}_{0}}] \}\nonumber\\
&+&|\mathbb{P}^{\epsilon}_{\kappa_{0},0}[(X^{\epsilon}_{\tau(-\kappa,\kappa^{'}_{0})},Y^{\epsilon}_{\tau(-\kappa,\kappa^{'}_{0})}) \in \Gamma_{\kappa^{'}_{0}}]-\mathbb{P}^{\epsilon}_{-\kappa_{0},0}[(X^{\epsilon}_{\tau(-\kappa,\kappa^{'}_{0})},Y^{\epsilon}_{\tau(-\kappa,\kappa^{'}_{0})}) \in \Gamma_{\kappa^{'}_{0}}]|\times\nonumber\\
& &\times\max\{\sup_{y\in D^{\epsilon}_{\kappa^{'}_{0}}} F^{\epsilon}(\kappa^{'}_{0},y),\inf_{y\in D^{\epsilon}_{\kappa^{'}_{0}}} F^{\epsilon}(\kappa^{'}_{0},y)\}\nonumber\\
&\leq&|\sup_{y\in D^{\epsilon}_{\kappa^{'}_{0}}} F^{\epsilon}(\kappa^{'}_{0},y)-\inf_{y\in D^{\epsilon}_{\kappa^{'}_{0}}} F^{\epsilon}(\kappa^{'}_{0},y)|+\nonumber\\
&+&|\mathbb{P}^{\epsilon}_{\kappa_{0},0}[(X^{\epsilon}_{\tau(-\kappa,\kappa^{'}_{0})},Y^{\epsilon}_{\tau(-\kappa,\kappa^{'}_{0})}) \in \Gamma_{\kappa^{'}_{0}}]-\mathbb{P}^{\epsilon}_{-\kappa_{0},0}[(X^{\epsilon}_{\tau(-\kappa,\kappa^{'}_{0})},Y^{\epsilon}_{\tau(-\kappa,\kappa^{'}_{0})}) \in \Gamma_{\kappa^{'}_{0}}]|\label{MaxMinProbAlternativeStatement2b1}
\end{eqnarray}
Relation (\ref{MaxMinProbAlternativeStatement1a}) holds with
$\kappa^{'}_{0}$ in place of $\kappa_{0}$ as well. Namely, for
$(\kappa^{'}_{0},y_{1}),(\kappa^{'}_{0},y_{2})\in \Gamma_{\kappa^{'}_{0}}$ we have
\begin{displaymath}
|F^{\epsilon}(\kappa^{'}_{0},y_{1})-F^{\epsilon}(\kappa^{'}_{0},y_{2})|\rightarrow
0 \textrm{ as } \epsilon\downarrow 0,
\end{displaymath}
uniformly in $y_{1},y_{2}\in D^{\epsilon}_{\kappa^{'}_{0}}$. This implies that the first term on the right hand side of (\ref{MaxMinProbAlternativeStatement2b1}) can be made arbitrarily small.

We show now how the second term on the right hand side of (\ref{MaxMinProbAlternativeStatement2b1}) can be made arbitrarily small.
For $\epsilon$ sufficiently small and for $\kappa^{'}_{0}$ much smaller than a small $\kappa$, it can be shown that $\mathbb{P}^{\epsilon}_{\kappa_{0},0}[(X^{\epsilon}_{\tau(-\kappa,\kappa^{'}_{0})},Y^{\epsilon}_{\tau(-\kappa,\kappa^{'}_{0})}) \in \Gamma_{\kappa^{'}_{0}}]$ and $\mathbb{P}^{\epsilon}_{-\kappa_{0},0}[(X^{\epsilon}_{\tau(-\kappa,\kappa^{'}_{0})},Y^{\epsilon}_{\tau(-\kappa,\kappa^{'}_{0})}) \in \Gamma_{\kappa^{'}_{0}}]$ are arbitrarily close to 1.  This can be done by an argument similar to the one that was used to prove Lemma \ref{Lemma24}. Similarly to there, one estimates the difference $$|\mathbb{P}^{\epsilon}_{\kappa_{0},0}[(X^{\epsilon}_{\tau(-\kappa,\kappa^{'}_{0})},Y^{\epsilon}_{\tau(-\kappa,\kappa^{'}_{0})}) \in \Gamma_{\kappa^{'}_{0}}]-\hat{\mathbb{P}}^{\epsilon}_{\kappa_{0}}[\hat{X}^{\epsilon}_{\hat{\tau}(-\kappa,\kappa^{'}_{0})}={\kappa^{'}_{0}}]|$$ and uses the corresponding estimate for $\hat{\mathbb{P}}^{\epsilon}_{\kappa_{0}}[\hat{X}^{\epsilon}_{\hat{\tau}(-\kappa,\kappa^{'}_{0})}={\kappa^{'}_{0}}]$ for $\epsilon$ sufficiently small, where $\hat{X}^{\epsilon}_{t}$ is the process defined by (\ref{EscordProcess}). One also needs to use Lemma \ref{Lemma24}. The treatment of $\mathbb{P}^{\epsilon}_{-\kappa_{0},0}[(X^{\epsilon}_{\tau(-\kappa,\kappa^{'}_{0})},Y^{\epsilon}_{\tau(-\kappa,\kappa^{'}_{0})}) \in \Gamma_{\kappa^{'}_{0}}]$ is almost identical with the obvious changes. We will not repeat the lengthy, but straightforward calculations here. Hence, the second term on the right hand side of (\ref{MaxMinProbAlternativeStatement2b1}) can be made arbitrarily small.

Using the above we finally get that
\begin{equation}
|F^{\epsilon}(\kappa_{0},0)-F^{\epsilon}(-\kappa_{0},0)| \rightarrow
0 \textrm{ as } \epsilon\downarrow 0.
\end{equation}
Therefore, we have established that (\ref{MaxMinProbAlternativeStatement}) holds for any $(x_{1},y_{1}),(x_{2},y_{2})\in \Gamma_{\pm\kappa_{0}}$.

Next, we prove that
\begin{equation}
\max_{(x,y)\in
\Omega_{0}}\mathbb{P}^{\epsilon}_{x,y}(X^{\epsilon}_{\tau(\pm
\kappa)} =\kappa)-\min_{(x,y)\in
\Omega_{0}}\mathbb{P}^{\epsilon}_{x,y}(X^{\epsilon}_{\tau(\pm
\kappa)}=\kappa)\rightarrow 0 \textrm{ as } \epsilon\downarrow
0,\label{MaxMinProb}
\end{equation}
where $\Omega_{0}=\lbrace(x,y):x\in[-\kappa_{0},\kappa_{0}], y\in
D^{\epsilon}_{x}\rbrace$. We use again the strong Markov property.
Let us choose some $(x,y)$ such that $|x|<\kappa_{0}$ and $y\in
D^{\epsilon}_{x}$. By strong Markov property with respect to the
exit time from $\{(x,y):|x|<\kappa_{0},y\in D^{\epsilon}_{x}\}$ we
have:
\begin{displaymath}
\inf_{(x,y)\in \Gamma_{\pm\kappa_{0}}}\mathbb{E}^{\epsilon}_{x,y}\chi_{(X^{\epsilon}_{\tau(\pm
\kappa)}=\kappa)}\leq \mathbb{P}^{\epsilon}_{x,y}(X^{\epsilon}_{\tau(\pm
\kappa)} =\kappa)\leq \sup_{(x,y)\in \Gamma_{\pm\kappa_{0}}}\mathbb{E}^{\epsilon}_{x,y}\chi_{(X^{\epsilon}_{\tau(\pm
\kappa)}=\kappa)}
\end{displaymath}
The latter implies that for any
$(x_{1},y_{1}),(x_{2},y_{2})\in\Omega_{0}$:
\begin{eqnarray}
|\mathbb{P}^{\epsilon}_{x_{1},y_{1}}(X^{\epsilon}_{\tau(\pm
\kappa)} =\kappa)&-&\mathbb{P}^{\epsilon}_{x_{2},y_{2}}(X^{\epsilon}_{\tau(\pm
\kappa)}=\kappa)|\leq \nonumber\\
&\leq& |\sup_{(x,y)\in \Gamma_{\pm\kappa_{0}}}\mathbb{E}^{\epsilon}_{x,y}\chi_{(X^{\epsilon}_{\tau(\pm
\kappa)}=\kappa)}
-\inf_{(x,y)\in \Gamma_{\pm\kappa_{0}}}\mathbb{E}^{\epsilon}_{x,y}\chi_{(X^{\epsilon}_{\tau(\pm
\kappa)})}|\nonumber
\end{eqnarray}
The latter inequality and (\ref{MaxMinProbAlternativeStatement})
imply (\ref{MaxMinProb}). Therefore
$\mathbb{P}^{\epsilon}_{x,y}(X^{\epsilon}_{\tau(\pm \kappa)}
=\kappa)$ has approximately the same value for all $(x,y)\in
\Omega_{0}$  for $\epsilon$ small enough. Let us now identify this
value. In order to do this we use the machinery with the invariant
measures. In particular, we consider again the two cycles of Markov
times $\lbrace \sigma_{n}\rbrace$ and $\lbrace \tau_{n}\rbrace$ that
are defined by (\ref{MarkovTimes3}).

The numbers $\nu_{\epsilon}^{\kappa_{0}}(\Gamma_{\pm\kappa_{0}})$ and $\nu_{\epsilon}^{\kappa}(\Gamma_{\pm\kappa})$ are strictly positive. This is because, starting from each of these sets, there is a positive probability  to reach every other $\Gamma_{\pm\kappa_{0}}$,$\Gamma_{\pm\kappa}$ in a finite number of cycles. We introduce the averages:
\begin{equation}
p_{+}(\epsilon)=\frac{\int_{\Delta_{\kappa_{0}}}\nu_{\epsilon}^{\kappa_{0}}(dx,dy)\mathbb{P}^{\epsilon}_{x,y}[(X^{\epsilon}_{\tau_{1}},Y^{\epsilon}_{\tau_{1}})\in
\Gamma_{\kappa}]}{\nu_{\epsilon}^{\kappa}(\Gamma_{\pm\kappa})}.\label{MaxProbability1}
\end{equation}
By relation (\ref{MaxMinProb}) we know that
\begin{equation}
|\mathbb{P}^{\epsilon}_{x,y}(X^{\epsilon}_{\tau(\pm
\kappa)} =\kappa)-p_{+}(\epsilon)|\leq \eta \label{MaxProbability2}
\end{equation}
for all $(x,y)$ such that $(x,y)\in [-\kappa_{0},  \kappa_{0}]\times
D^{\epsilon}_{x}$  for $\epsilon$ sufficiently small.

Moreover using the first of equations (\ref{InvariantMeasuresForMarkovChains}) we see that (\ref{MaxProbability1}) can be written as
\begin{equation}
p_{+}(\epsilon)=\frac{\nu_{\epsilon}^{\kappa}(\Gamma_{\kappa})}{\nu_{\epsilon}^{\kappa}(\Gamma_{\pm\kappa})}.\label{MaxProbability3}
\end{equation}
Furthermore, we for $0<\kappa_{0}<\kappa$ sufficiently small we have
\begin{displaymath}
u(\kappa)-u(\kappa_{0})\sim u'(0+)(\kappa-\kappa_{0}).
\end{displaymath}
The latter and Lemma \ref{Lemma52} imply for sufficiently small $0<\kappa_{0}<\kappa$ and sufficiently small $\epsilon$ that
\begin{equation}
\nu_{\epsilon}^{\kappa}(\Gamma_{\kappa})\sim\frac{1}{u'(0+)}
\frac{\epsilon}{\kappa-\kappa_{0}}.\label{AsymptoticsOfInvariantMeasures3}
\end{equation}
Similarly, we have for sufficiently small $0<\kappa_{0}<\kappa$ and sufficiently small $\epsilon$ that
\begin{equation}
\nu_{\epsilon}^{\kappa}(\Gamma_{-\kappa})\sim\frac{1}{u'(0-)}
\frac{\epsilon}{\kappa-\kappa_{0}}.\label{AsymptoticsOfInvariantMeasures3a}
\end{equation}
Therefore, we finally get for sufficiently small $0<\kappa_{0}<\kappa$ and sufficiently small $\epsilon$ that
\begin{equation}
\nu_{\epsilon}^{\kappa}(\Gamma_{\pm\kappa})\sim[\frac{1}{u'(0+)}+\frac{1}{u'(0-)}]
\frac{\epsilon}{\kappa-\kappa_{0}}.\label{AsymptoticsOfInvariantMeasures3b}
\end{equation}
Hence, by equations (\ref{MaxProbability3})-(\ref{AsymptoticsOfInvariantMeasures3b}) we get for sufficiently small $0<\kappa_{0}<\kappa$ and sufficiently small $\epsilon$ that
\begin{displaymath}
|p_{+}(\epsilon)-\frac{[u'(0+)]^{-1}}{[u'(0+)]^{-1}+[u'(0-)]^{-1}}|<\eta
\end{displaymath}
Similarly, we can prove that
\begin{displaymath}
 \max_{(x,y)\in
\Omega_{0}}\mathbb{P}^{\epsilon}_{x,y}(X^{\epsilon}_{\tau(\pm
\kappa)} =-\kappa)-\min_{(x,y)\in
\Omega_{0}}\mathbb{P}^{\epsilon}_{x,y}(X^{\epsilon}_{\tau(\pm
\kappa)}=-\kappa)\rightarrow 0 \textrm{ as } \epsilon\downarrow
0.
\end{displaymath}
Then for
\begin{displaymath}
p_{-}(\epsilon)=\frac{\nu_{\epsilon}^{\kappa}(\Gamma_{-\kappa})}{\nu_{\epsilon}^{\kappa}(\Gamma_{\pm\kappa})},
\end{displaymath}
we can obtain that
\begin{displaymath}
|p_{-}(\epsilon)-\frac{[u'(0-)]^{-1}}{[u'(0+)]^{-1}+[u'(0-)]^{-1}}|<\eta
\end{displaymath}
and that the corresponding relation (\ref{MaxProbability2}) holds, i.e.
\begin{displaymath}
|\mathbb{P}^{\epsilon}_{x,y}(X^{\epsilon}_{\tau(\pm
\kappa)} =-\kappa)-p_{-}(\epsilon)|\leq \eta
\end{displaymath}
for all $(x,y)$ such that $(x,y)\in [-\kappa_{0},  \kappa_{0}]\times
D^{\epsilon}_{x}$  for $\epsilon$ sufficiently small.

The latter concludes the proof of the Lemma.
\end{proof}
\section{Acknowledgements}
I would like to thank Professor Mark Freidlin for posing the
problem and for helpful discussions and Professor Sandra Cerrai, Anastasia Voulgaraki and Hyejin Kim for helpful
suggestions. I would also like to thank the anonymous referee for the constructive comments and suggestions that greatly improved the paper.

\end{document}